\documentclass{article}

\pdfoutput=1

\usepackage{graphicx}
\usepackage{amsmath}
\usepackage{amssymb}
\usepackage{amsthm}
\usepackage[utf8]{inputenc}
\usepackage{hyperref}
\hypersetup{
    pdftitle={Long-time existence of Brownian motion on configurations of two landmarks},
    pdfauthor={Karen Habermann, Philipp Harms, Stefan Sommer},
}

\newcommand{\Diff}{\mathrm{Diff}}
\newcommand{\Id}{\operatorname{Id}}

\newcommand{\R}{\mathbb{R}}
\newcommand{\N}{\mathbb{N}}

\newcommand{\E}{\mathbb{E}}

\newcommand{\e}{\operatorname{e}}
\newcommand{\dd}{\,{\mathrm d}}
\newcommand{\db}{{\mathrm d}}
\newcommand{\pt}{\partial}

\newcommand{\cc}{\lambda}

\newtheorem*{theorem}{Theorem}

\title{Long-time existence of Brownian motion on configurations of two landmarks\thanks{The authors thank Marc Arnaudon for helpful discussions and comments. Philipp Harms gratefully acknowledges financial support by the National Research Foundation Singapore under the award NRF-NRFF13-2021-0012 and by Nanyang Technological University Singapore under the award NAP-SUG. Stefan Sommer is supported by the Villum Foundation grant 40582, the Novo Nordisk Foundation grant NNF18OC0052000, and the UCPH Data+ strategy funds.}}
\author{Karen Habermann, Philipp Harms, Stefan Sommer}
\date{\today}

\begin{document}

\maketitle

\begin{abstract}
    We study Brownian motion on the space of distinct landmarks in $\R^d$, considered as a homogeneous space with a Riemannian metric inherited from a right-invariant metric on the diffeomorphism group. As of yet, there is no proof of long-time existence of this process, despite its fundamental importance in statistical shape analysis, where it is used to model stochastic shape evolutions. We make some first progress in this direction by providing a full classification of long-time existence for configurations of exactly two landmarks, governed by a radial kernel. For low-order Sobolev kernels, we show that the landmarks collide with positive probability in finite time, whilst for higher-order Sobolev and Gaussian kernels, the landmark Brownian motion exists for all times. We illustrate our theoretical results by numerical simulations. 
\end{abstract}

\section{Introduction}

A common approach to shape analysis is to model shape variations as diffeomorphic deformations of the ambient space, where the shapes reside. Then, shape variations can be quantified using a right-invariant Riemannian metric on a diffeomorphism group, see Younes~\cite{younes}. In particular, one obtains in this way Riemannian metrics on landmark spaces which are well described by a kernel, as discussed by Micheli in~\cite{micheliDifferentialGeometryLandmark2008}. The metric is geodesically or metrically complete if the associated reproducing kernel Hilbert space embeds continuously in $C^1_b$ or $C^2_b$, respectively, see Bauer, Bruveris and Michor~\cite{bauerOverviewGeometriesShape2014}. Some weaker sufficient conditions are known as well, for example, as in Joshi and Miller~\cite{joshi2000landmark}. 

Brownian completeness is a related but distinct notion, which concerns the long-time existence or equivalently the non-explosion of Brownian motion on a Riemannian manifold. For a general overview, one may consult the monographs by Hackenbroch and Thalmaier~\cite{hackenbroch1994stochastische} or Hsu~\cite{hsu}. 
Brownian completeness is of fundamental importance in statistical shape analysis, where Brownian motion is often the first stochastic process to consider when modelling stochastic shape evolutions. For instance, it has been used for statistics of shapes in \cite{sommerBridgeSimulationMetric2017, stanevaLearningShapeTrends2017} and for diffusion means of shapes in \cite{eltznerDiffusionMeansGeometric2022}. Brownian completeness is (implicitly) assumed in these works but is unknown even in the simplest case of landmark spaces. As a consequence, these statistical works are built on an uncertain foundation, and establishing non-explosion criteria for Riemannian Brownian motion on landmark spaces has direct applied impact. 

In this article, we fully characterise Brownian completeness in the case of two landmarks in $\R^d$, for $d\geq 1$. In particular, we establish the following result, stated in terms of notions and notations rigorously introduced in Section~\ref{sec:lddmm} and Section~\ref{sec:kernels}.
\begin{theorem}
    Let $Q$ be the landmark manifold of pairs of distinct points in $\R^d$, for $d\geq 1$. Let $g$ be the Riemannian metric on $Q$, whose inverse is given by 
    \begin{displaymath}
        g^{-1}_q(\xi,\eta)
        =
        \sum_{i,j=1}^2 \xi_i^\top K(q_i,q_j)\eta_j\;,
        \qquad
        q \in Q\;, 
        \quad 
        \xi, \eta \in T^*_qQ\;,
    \end{displaymath}
    for some positive definite radial kernel 
    \begin{displaymath}
        K\colon \R^d\times\R^d \to \R^{d\times d},
        \qquad
        K(x,y)=k(\|x-y\|_{\R^d}) I_d\;,
    \end{displaymath}
    where $\|\cdot\|_{\R^d}$ is the Euclidean norm on $\R^d$,
    $I_d$ is the $d\times d$ identity matrix, and $k\colon(0,\infty)\to\R$. 
    Suppose that $k$ extends continuously to $[0,\infty]$, vanishes at $\infty$, is continuously differentiable on $(0,\infty)$, and has a bounded and Lipschitz continuous derivative on $[1,\infty)$. 
    Moreover, suppose that, for $D, \gamma >0$,
    \begin{displaymath}
        k(0)-k(r)=Dr^\gamma +o\left(r^\gamma\right)\;,
        \qquad
        \text{as $r\downarrow 0$}\;.
    \end{displaymath}
    Then, the Riemannian manifold $(Q,g)$ is Brownian complete if $\gamma\geq 2$, whilst it is Brownian incomplete if $\gamma<2$.
\end{theorem}

The theorem is proven in Section~\ref{sec:twolandmarks}.
It implies that low-order Sobolev metrics on the space of two landmarks are Brownian incomplete, whereas high-order Sobolev metrics, and also the metric with Gaussian kernel, are Brownian complete. 
Beyond characterising Brownian completeness, our analysis provides a detailed description of the long-term behaviour of the two-landmark system. 
These results are presented in Section~\ref{sec:application}.

It is important to remark that these results have no direct implications to or from metric completeness. Indeed, there are Brownian complete yet metrically incomplete spaces such as the punctured plane, and conversely, there are also Brownian incomplete yet metrically complete spaces, see e.g.~\cite{hackenbroch1994stochastische}. Nevertheless, the above examples show that higher-order metrics have favourable Brownian completeness properties, which is in line with similar results for metric completeness. 

To summarise, the present article takes some first steps towards solving the general question of Brownian completeness of landmark spaces. Whilst our analysis is presently limited to configurations of merely two landmarks, we do hope that our ideas and techniques will eventually lead to a solution of the Brownian completeness question for more general shape spaces.

\subsection*{Overview of the paper}
The article is organised as follows. Section~\ref{sec:landmarkBM} sets the stage by introducing the Riemannian landmark manifold, the Riemannian Brownian motion thereon, and several radial kernels of interest. This is followed in Section~\ref{sec:twolandmarks} by the analysis of the long-time existence of Riemannian Brownian motion on the space of landmark pairs. In Section~\ref{sec:distantceprocess}, we derive the It\^o stochastic differential equation for the inter-landmark distance process. In Section~\ref{sec:collision}, we analyse the singularity of the distance process at zero and characterise the possibility of landmark collision in finite time. Section~\ref{sec:escape} rules out the possibility of landmarks escaping to infinity before collision. Section~\ref{sec:application} summarises the preceding analysis and provides a fine-grained description of the inter-landmark distance process. Finally, in Section~\ref{sec:numerics}, we show that our theoretical results are in line with numerical simulations.

\section{Landmark Brownian motion}\label{sec:landmarkBM}

\subsection{Landmark space}\label{sec:lddmm}

In the Large Deformation Diffeomorphic Metric Mapping (LDDMM) framework, see~\cite{younes}, shape variations are modelled as diffeomorphic deformations. This makes the framework applicable to a wide range of shape spaces. An important shape space within this framework is the landmark manifold $Q$, which consists of configurations of $n\geq 2$ distinct landmark points in $\mathbb R^{d}$, for $d\geq 1$. Thus, $Q$ could be seen as an open subset of $\R^{nd}$. In contrast, the LDDMM framework views it as a Riemannian homogeneous space, whose metric stems from a right-invariant metric on a group of diffeomorphisms on $\R^d$. This perspective originates from diffeomorphic landmark matching, where two landmark configurations $q,p\in Q$ are matched by solving the following optimisation problem:
\begin{equation}
    \operatorname*{minimise}_u
    \int_0^1 \|u_t\|_V^2 \dd t
    \quad\operatorname{subject\ to}\quad
    \varphi_1.q=p\;.
        \label{eq:LDDMM_energy}
\end{equation} 
Here, $u\colon [0,1]\to \mathfrak X_c(\mathbb R^d)$ is a  time-dependent compactly supported vector field which generates a flow $\varphi\colon [0,1]\to\Diff_c(\mathbb R^d)$ of compactly supported diffeomorphisms via the flow equation
\begin{equation}
    \frac{\partial\varphi_t}{\partial t}
    =
    u_t\circ\varphi _t\;,
    \qquad
    \varphi_0=\Id_{\mathbb R^d}\;.
    \label{eq:LDDMM_rec}
\end{equation}
Instead of compact support, one may impose alternative regularity conditions such as rapid decay, boundedness of all derivatives, or quasi-analyticity, see~\cite{michor2013zoo} and~\cite{kriegl2015exotic}. 
The end point of the diffeomorphic flow is denoted by $\varphi_1$, and the constraint in \eqref{eq:LDDMM_energy} forces $\varphi_1.q:=(\varphi_1(q_1),\dots,\varphi_1(q_n))$ to match up with $p$.
The choice of norm $\|\cdot\|_V$ is discussed in detail in Section~\ref{sec:kernels}. 
For now, it suffices to assume that it stems from an inner product $\langle\cdot,\cdot\rangle_V$ on the space $\mathfrak X_c(\R^d)$ such that point evaluations are norm-continuous. 
Then, the completion of $\mathfrak X_c(\R^d)$ with respect to this norm is a Hilbert space $V$ with positive reproducing kernel $K\colon \R^d\times\R^d\to\R^{d\times d}$. The kernel gives rise to an integral operator $K\colon V^*\to V$, whose inverse is the inner product $\langle\cdot,\cdot\rangle_V\colon V\to V^*$. See \cite{younes} for further details.

A geometric interpretation is as follows. The space $\mathfrak X_c(\R^d)$ of compactly supported vector fields is the Lie algebra of the infinite-dimensional Lie group $\Diff_c(\R^d)$ of compactly supported diffeomorphisms, as discussed by Michor and Mumford~\cite{michor2013zoo}. The inner product $\langle\cdot,\cdot\rangle_V$ on $\mathfrak X_c(\R^d)$ can be extended to a unique right-invariant weak Riemannian metric $g$ on the Lie group $\Diff_c(\R^d)$. Compactly supported diffeomorphisms $\varphi\in\Diff_c(\R^d)$ act on landmark configurations $q =(q_1,\dots,q_n) \in Q$ from the left as
\begin{displaymath}
    \varphi.q=(\varphi(q_1),\ldots,\varphi(q_n))\;.
\end{displaymath}
For fixed $q$ and variable $\varphi$, this action is a submersion from $\Diff_c(\R^d)$ to $Q$. 
There is a unique Riemannian metric $g$ on $Q$ such that this submersion is Riemannian, see~\cite{micheli2012sectional}. 
Moreover, the optimisation problem \eqref{eq:LDDMM_energy} is equivalent to the geodesic boundary value problem for this metric, that is, the infimal energy in \eqref{eq:LDDMM_energy} is the squared geodesic distance between $q$ and $p$, and any minimiser $u$ generates a diffeomorphic flow $\varphi$ which projects down to a geodesic in $Q$.
If $V$ embeds continuously in $\mathfrak X_{C^1_b}(\R^d)$, then $\Diff_c(\R^d)$ can be completed to the half-Lie group $\Diff_V(\R^d)$, see~\cite{bauer2023regularity}, which is modeled on $V$, carries the strong Riemannian metric extended right-invariantly from $\langle\cdot,\cdot\rangle_V$, and is geodesically complete, see~\cite{younes}. Consequently, $Q$ is geodesically complete. 
If $V$ embeds continuously in $\mathfrak X_{C^2_b}(\R^d)$, then the metric on $Q$ is $C^2$, and the Hopf--Rinow theorem implies that $Q$ is also metrically complete, see~\cite{micheli2012sectional}. Further references for these arguments are provided in the overview article \cite{bauerOverviewGeometriesShape2014}.

Computationally, it is important that the cometric $g^{-1}$, that is, the inverse of the Riemannian metric $g$, admits a simple description in terms of the reproducing kernel $K$, see~\cite{micheli2012sectional}, namely, for $q\in Q$ and covectors $\xi,\eta\in T_q^*Q$,
\begin{equation}\label{eq:cometric}
    g_q^{-1}(\xi,\eta)
    =
    \sum_{i,j=1}^n\xi_i^\top K(q_i,q_j)\eta_j\;.
\end{equation}
Our analysis is based on this formula alone and does not make use of its geometric origins.
A concise derivation of this formula and elegant expressions for the corresponding geodesic equation and curvature can be found in \cite{michor2020manifolds}.

\subsection{Brownian motion of landmark configurations}

The Riemannian metric $g$ on the landmark configuration space $Q$ gives rise to the Laplace--Beltrami operator $\Delta_Q$ on $Q$. Brownian motion on $Q$ is the diffusion process $(q_t)_{t\in[0,\zeta)}$ on $Q$ generated by $\frac{1}{2}\Delta_Q$ with some initial value $q_0$ and defined for $t \geq 0$ up to some explosion time $\zeta\in(0,\infty]$, see~\cite{hackenbroch1994stochastische,hsu} for further details. Equivalently, as discussed by Hsu~\cite[Example 3.3.5]{hsu}, this diffusion process is the unique strong solution to a certain It\^o stochastic differential equation written in charts. On the landmark space $Q$, which is a subset of $\mathbb R^{nd}$, a single chart suffices. Hence, Brownian motion $(q_t)_{t\in[0,\zeta)}$ can be expressed as the unique strong solution to the It\^o stochastic differential equation, for $i\in\{1,\dots,n\}$,
\begin{equation}\label{eq:RiemBM}
    \db q_t^i
    =-\frac{1}{2} \sum_{\ell,m=1}^n K(q_t)^{\ell m}\Gamma(q_t)_{\ell m}^i \dd t
    + \sqrt{K(q_t)}^i\dd W_t
\end{equation}
subject to initial condition $q_0$ and defined for $t \in [0, \zeta)$, where $(W_t)_{t\geq 0}$ is an $\R^{nd}$-valued Wiener process on some stochastic basis $(\Omega,\mathcal F,(\mathcal F_t)_{t\geq 0},\mathbb P)$ satisfying the usual conditions.
Here, $K(q_t)^{\ell m}=K(q_t^\ell,q_t^m)$ denotes the cometric \eqref{eq:cometric} at $q_t \in Q$, with indices $\ell,m\in \{1,\dots,n\}$. Moreover, $\Gamma$ denotes the Christoffel symbol associated with the metric $g$.

According to~\cite{hackenbroch1994stochastische}, the Riemannian manifold $Q$ is called Brownian complete if, for every initial value $q_0\in Q$, the corresponding explosion time $\zeta=\zeta(q_0)$ satisfies
\begin{displaymath}
    \mathbb P(\zeta(q_0)=\infty)=1\;.
\end{displaymath}
Letting $p$ denote the Dirichlet heat kernel on $Q$ with respect to the Riemannian volume measure $\operatorname{vol}_g$, Brownian completeness is equivalent to
\begin{displaymath}
    \int_Q p(t,q_0,q)\operatorname{vol}_g(\db q) =1
    \quad\text{for all}\quad(t,q_0)\in (0,\infty)\times Q\;.
\end{displaymath}
There exist several sufficient conditions which guarantee that a Riemannian manifold is Brownian complete, such as the manifold being compact or the manifold having Ricci curvature bounded from below. For further details, one may consult~\cite{hackenbroch1994stochastische, hsu}. 
These conditions are difficult to apply in our setting because the landmark space $Q$ is non-compact, and our numerical simulations suggest that its Ricci curvature is unbounded from below. We therefore proceed with an alternative analysis, which results in a refined description of the long-term behaviour of Riemannian Brownian motion on landmark space, but is presently limited to configurations of two landmarks.

\subsection{Kernels}\label{sec:kernels}

We restrict our attention to kernels which are invariant under rotations and translations. This assumption is satisfied in the most important examples and significantly simplifies our analysis. Thus, we consider positive definite kernels of the form 
\begin{equation}\label{eq:radialkernel}
    K\colon \R^d\times\R^d \to\R^{d\times d}, 
    \qquad
    (q_i,q_j)\mapsto k(\|q_i-q_j\|_{\mathbb R^d}) I_d\;,
\end{equation}
where $I_d$ is the $d\times d$ identity matrix and $k\colon (0,\infty)\to\R$ is a scalar function.

An important special case are Bessel potentials of order $\alpha>d$, also known as Sobolev kernels, see Aronszajn and Smith~\cite{aronszajn1961theory}. 
The Bessel potential of order $\alpha \in (0,\infty)$ in $d\in\mathbb N$ dimensions is defined as the radial kernel \eqref{eq:radialkernel} with  
\begin{equation}\label{eq:bessel}
k(r) 
=
\frac{1}{2^{(\alpha-2)/2}(2\pi)^{d/2}\Gamma(\alpha/2)}r^\nu J_{-\nu}(r)\;,
\qquad
\nu=\frac{\alpha-d}{2}\;,
\end{equation}
where $J_{-\nu}=J_{\nu}$ on the right-hand side is the modified Bessel function, which is denoted by $K_{-\nu}=K_{\nu}$ in much of the literature. 
The integral operator with kernel $K$ is the inverse of the fractional Laplacian $$(\Id-\Delta)^{\alpha/2}\colon\mathfrak X_{H^{\alpha/2}}(\R^d)\to \mathfrak X_{L^2}(\R^d)\;,$$ as can be seen from \cite[Equation~(4, 6)]{aronszajn1961theory} in the Fourier domain. 
Moreover, for $\alpha>d$, the kernel $K$ is positive and it is the reproducing kernel of the Sobolev space $\mathfrak X_{H^{\alpha/2}}(\R^d)$, as established in~\cite{aronszajn1961theory}.
The asymptotics for small $r$ are given by Abramowitz and Stegun~\cite[Equations (6.1.15-17) and (9.6.2--11)]{abramowitz1972handbook} as
\begin{equation}\label{eq:bessel-asymptotics}
r^\nu J_{-\nu}(r)
=
2^{\nu-1}\Gamma(\nu)-
\begin{cases}
-2^{\nu-1}\Gamma(-\nu)r^{2\nu}+o\left(r^{2\nu}\right),
&\nu\in(0,1)\;,\\
-2^{-1}r^2\log(r)+o\left(r^2\log(r)\right),
&\nu=1\;,\\
2^{\nu-3}\Gamma(\nu-1)r^2+o\left(r^2\right),
&\nu\in (1,\infty)\;.
\end{cases}
\end{equation}
Whenever $\nu-1/2$ is a natural number, the kernel admits the explicit formula \cite[Lemma~9.16]{younes}, which gives rise to 
\begin{displaymath}
    k(r)
    \propto 
    \e^{-r} \sum_{l=0}^{\nu-1/2} \frac{2^l(2c-l)!}{(c-l)!\,l!} \,r^l\;,
\end{displaymath}
where $\propto$ denotes equality up to a positive constant.
Some special cases are 
\begin{equation*}
k(r) \propto \e^{-r} 
\begin{cases}
1\;, & \nu=\tfrac12\;, \\
1+r\;, & \nu=\tfrac32\;, \\
3+3r+r^2\;, & \nu=\tfrac52\;,\\
15+15r+6r^2+r^3\;, & \nu=\tfrac72\;.
\end{cases}
\end{equation*}

The Gaussian kernel, that is, the radial kernel~\eqref{eq:radialkernel} with $k(r)=\exp(-r^2)$
is another important example. It is positive definite and can be seen as a Bessel kernel of infinite order because its Fourier multiplier $\e^{-\pi^2\|\xi\|^2}$ is the limit $\alpha\to\infty$ of the rescaled Bessel Fourier multiplier $(1+\pi^2/\alpha\|\xi\|^2)^{-\alpha}$. 
The asymptotics for the Gaussian kernel for small $r$ are
\begin{displaymath}
    k(r)= 1-r^2 + o\left(r^3\right)\;.
\end{displaymath}

Interestingly, as shown next, different near-zero asymptotics of the kernel result in qualitatively different long-term behaviour for the Riemannian Brownian motion on landmark space. This analysis requires that $k$ extends continuously to $[0,\infty]$, vanishes at $\infty$, is continuously differentiable on $(0,\infty)$, and has a bounded and Lipschitz continuous derivative on $[1,\infty)$. 
All of these assumptions are satisfied for Bessel potentials of order $\alpha>d$ by Abramowitz and Stegun~\cite[Equations (9.6.1), (9.6.9), and (9.7.2--4)]{abramowitz1972handbook}, and are also satisfied for the Gaussian kernel.

\section{Brownian motion of two landmarks}\label{sec:twolandmarks}

We next characterise the long-term behaviour of Riemannian Brownian motion on configurations of exactly two landmarks, endowed with a Riemannian structure induced by a radial kernel as described in the preceding section. 
The key observation is that for a radial kernel, the distance between the two landmarks is a diffusion process, whose dynamics is characterised by a scalar stochastic differential equation. It then remains to study the singularity at zero of this one-dimensional diffusion process. For this, we follow the classification of singular points by Cherny and Engelbert in~\cite{SSDE}.

To derive the It\^o stochastic differential equation for the distance process, we significantly reduce the complexity of the required computations by working in a well chosen coordinate system. Specifically, unlike a brute-force application of It\^o's formula to the stochastic differential equation \eqref{eq:RiemBM}, our approach circumvents the need to determine all Christoffel symbols. Instead, it requires only the computation of one divergence. The reduction in complexity is best illustrated in case of two landmarks in $\R$, which is why we first discuss this case, even though the result is included in the subsequent more general analysis for two landmarks in $\R^d$.

\subsection{Distance process between two landmarks}\label{sec:distantceprocess}

Let $K\colon \R^d\times\R^d\to\R^{d\times d}$ be a radial kernel~\eqref{eq:radialkernel} described in terms of a functions $k\colon [0,\infty)\to\R$ which is continuous on $[0,\infty)$ and continuously differentiable on $(0,\infty)$. We further set $\cc=k(0)$.

\subsubsection{Configurations with two landmarks in \texorpdfstring{$\R$}{R}}
For landmark configurations with two landmarks in $\R$, we have
\begin{displaymath}
    Q=\{q=(x,y):x,y\in\R\text{ such that }x\not=y\}\;.
\end{displaymath}
Since the radial kernel $K$ takes the form~\eqref{eq:radialkernel}, we have, for $q=(x,y)\in Q$,
\begin{displaymath}
    K(q)=
    \begin{pmatrix}
        \cc & k(|x-y|)\\
        k(|x-y|) & \cc
    \end{pmatrix}\;.
\end{displaymath}
Due to~\eqref{eq:cometric} describing the cometric induced by the Green's kernel $K$, it follows that the metric $g$ on $Q$ induced by $K$ is given as
\begin{displaymath}
    g_q=K(q)^{-1}
    =\frac{1}{\cc^2-(k(|x-y|))^2}
    \begin{pmatrix}
        \cc & -k(|x-y|) \\
        -k(|x-y|) & \cc
    \end{pmatrix}\;.
\end{displaymath}
We now change to a system of coordinates in which the expression for the metric $g$ diagonalises. For $q=(x,y)\in Q$, we set
\begin{displaymath}
    u=x-y
    \quad\text{and}\quad
    v=x+y\;.
\end{displaymath}
The constraint $x\not=y$ then amounts to the condition $u\not= 0$. Without loss of generality, we may work in the half plane $u>0$, which corresponds to the assumption that the Riemannian Brownian motion is started from a landmark configuration where the first landmark is bigger than the second one. The Brownian motion on $(Q,g)$ only leaves the half plane defined by $u>0$ if the two landmarks collide in finite time, that is, the Riemannian Brownian motion explodes.

From $2x=u+v$ and $2y=v-u$, we obtain
\begin{displaymath}
  \frac{\pt}{\pt u}=
  \frac{1}{2}\left(\frac{\pt}{\pt x}-\frac{\pt}{\pt y}\right)
  \quad\text{and}\quad
  \frac{\pt}{\pt v}=
  \frac{1}{2}\left(\frac{\pt}{\pt x}+\frac{\pt}{\pt y}\right)\;.
\end{displaymath}
Since $|x-y|=u$ in the half plane $u>0$, it follows that, for $q=(u,v)\in Q$,
\begin{equation}\label{eq:uvnorms}
    g_q\left(\frac{\pt}{\pt u},\frac{\pt}{\pt u}\right)
    =\frac{1}{2\left(\cc-k(u)\right)}
    \quad\text{and}\quad
    g_q\left(\frac{\pt}{\pt v},\frac{\pt}{\pt v}\right)
    =\frac{1}{2\left(\cc+k(u)\right)}
\end{equation}
as well as
\begin{displaymath}
    g_q\left(\frac{\pt}{\pt u},\frac{\pt}{\pt v}\right)=0\;.
\end{displaymath}
Therefore, the vector fields $X_1$ and $X_2$ on $Q$ defined by
\begin{equation}\label{eq:onframe}
  X_1=\sqrt{2(\cc-k(u))}\,\frac{\pt}{\pt u}
  \quad\text{and}\quad
  X_2=\sqrt{2(\cc+k(u))}\,\frac{\pt}{\pt v}  
\end{equation}
form an orthonormal frame $(X_1,X_2)$ for the tangent bundle $TQ$ with respect to the metric $g$. In particular, we can write
\begin{equation}\label{eq:niceLB}
    \Delta_Q=X_1^2+X_2^2+(\operatorname{div} X_1)X_1 +
  (\operatorname{div} X_2)X_2\;,
\end{equation}
where the divergence is taken with respect to the induced Riemannian volume measure. From the expression~\eqref{eq:niceLB}, we can read off that the Brownian motion on $(Q,g)$, that is, the diffusion process $(q_t)_{t\in[0,\zeta)}$ on $Q$ with generator $\frac{1}{2}\Delta_Q$, is the unique strong solution to the Stratonovich stochastic differential equation
\begin{displaymath}
    \db q_t=X_1(q_t)\circ\dd B_t + X_2(q_t)\circ\dd W_t +
    \frac{1}{2}(\operatorname{div} X_1)X_1(q_t)\dd t +
    \frac{1}{2}(\operatorname{div} X_2)X_2(q_t)\dd t\;,
\end{displaymath}
where $(B_t)_{t\geq 0}$ and $(W_t)_{t\geq 0}$ are independent one-dimensional standard Brownian motions.

Due to the form~\eqref{eq:onframe} of the vector field $X_1$, it follows that, for $q_t=(u_t,v_t)$, the distance process $(u_t)_{t\in[0,\zeta)}$ between the two landmarks induced by the Brownian motion on $(Q,g)$ is the unique strong solution to the Stratonovich stochastic differential equation
\begin{equation}\label{eq:2landStrat}
    \db u_t=\sqrt{2(\cc-k(u_t))}\circ\dd B_t +
    \frac{1}{2}(\operatorname{div} X_1)(u_t)\sqrt{2(\cc-k(u_t))}\dd t\;.
\end{equation}
It remains to compute the divergence of the vector field $X_1$ explicitly and to express the Stratonovich stochastic differential equation as an It\^o stochastic differential equation.

As a consequence of~\eqref{eq:uvnorms}, the Riemannian volume form $\operatorname{vol}_g$ on $Q$ induced by the Riemannian metric $g$ can be expressed in the coordinates $(u,v)$ as
\begin{displaymath}
    \operatorname{vol}_g=\frac{1}{4(\cc^2-(k(u))^2)}\dd u\dd v\;.
\end{displaymath}
It follows that
\begin{displaymath}
    (\operatorname{div} X_1)(u)
    =2\sqrt{\cc^2-(k(u))^2}\,
    \frac{\pt}{\pt u}\left(\frac{1}{\sqrt{2(\cc+k(u))}}\right)
    =-\frac{k'(u)\sqrt{\cc-k(u)}}{\sqrt{2}\,(\cc+k(u))}\;,
\end{displaymath}
which yields
\begin{displaymath}
    (\operatorname{div} X_1)(u)\sqrt{2(\cc-k(u))}
    =-\frac{k'(u)(\cc-k(u))}{\cc+k(u)}\;.
\end{displaymath}
We further note that, for $f\in C^2(Q)$,
\begin{displaymath}
    X_1(X_1(f))=2(\cc-k)\frac{\pt^2 f}{\pt u^2} - \frac{\pt k}{\pt u}\frac{\pt f}{\pt u}\;,
\end{displaymath}
from which we can read off the drift term contribution arising from $X_1^2$.
Since
\begin{displaymath}
    \frac{1}{2}\left((\operatorname{div} X_1)(u)\sqrt{2(\cc-k(u))}-k'(u)\right)
    =-\frac{\cc\, k'(u)}{\cc+k(u)}\;,
\end{displaymath}
we deduce that the Stratonovich stochastic differential equation~\eqref{eq:2landStrat} can be rewritten as the It\^o stochastic differential equation
\begin{equation}\label{eq:distin1D}
    \db u_t=\sqrt{2(\cc-k(u_t))}\dd B_t - \frac{\cc\, k'(u_t)}{\cc+k(u_t)}\dd t\;.
\end{equation}
One can check that this is consistent with the expression obtained by starting directly from~\eqref{eq:RiemBM}.

\subsubsection{Configurations with two landmarks in \texorpdfstring{$\R^d$}{Rd}}

The restriction to radial kernels of the form~\eqref{eq:radialkernel} allows us to describe the distance process between two landmarks, provided no additional landmarks are present, by a one-dimensional stochastic differential equation. This argument, which was developed for two landmarks in $\R$ in the previous subsection, carries over to two landmarks in $\R^d$, as shown next. 

As before, we simplify the computations significantly by working in suitable coordinates, for which the metric tensor diagonalises. When considering configurations consisting of two landmarks in $\R^d$, for $d\geq 1$, we have
\begin{displaymath}
    Q=\{q=(x,y):x,y\in\R^d\text{ such that }x\not=y\}\;.
\end{displaymath}
Moreover, the radial kernel $K$ of the form~\eqref{eq:radialkernel} is given by, for $q=(x,y)\in Q$,
\begin{displaymath}
    K(q)=
    \begin{pmatrix}
        \cc I_d & k(\|x-y\|_{\R^d}) I_d\\
        k(\|x-y\|_{\R^d}) I_d & \cc I_d
    \end{pmatrix}\;,
\end{displaymath}
where $I_d$ denotes the $d\times d$ identity matrix. As in the analysis for two landmarks in $\R$, it further follows from~\eqref{eq:cometric} and the above expression for the radial kernel $K$ that the induced metric $g$ on $Q$ is determined by
\begin{displaymath}
    g_q
    =K(q)^{-1}
    =\frac{1}{\cc^2-(k(\|x-y\|_{\R^d}))^2}
    \begin{pmatrix} 
        \cc I_d & -k(\|x-y\|_{\R^d})I_d \\
        -k(\|x-y\|_{\R^d}) I_d & \cc I_d
    \end{pmatrix}\;.
\end{displaymath}

We proceed by changing coordinates for the landmark space $Q$ from $(x,y)$ to $(u,v)$, where $u,v\in\R^d$ are given by
\begin{displaymath}
    u=x-y
    \quad\text{and}\quad
    v=x+y\;,
\end{displaymath}
and by setting
\begin{displaymath}
    r=\|x-y\|_{\R^d}=\|u\|_{\R^d}=\sqrt{\sum_{i=1}^d\left(u^i\right)^2}\;.
\end{displaymath}
Note that the constraint $x\not=y$ is equivalent to the condition $r\not =0$. As two landmarks collide if and only if their distance process hits zero, it suffices to study the stochastic dynamics of the distance process, for landmarks evolving according to a Riemannian Brownian motion. For $i,j\in\{1,\dots,d\}$, we obtain
\begin{displaymath}
    g_q\left(\frac{\pt}{\pt u^i},\frac{\pt}{\pt u^j}\right)
    =\frac{\delta_{ij}}{2\left(\cc-k(r)\right)}
    \quad\text{and}\quad
    g_q\left(\frac{\pt}{\pt v^i},\frac{\pt}{\pt v^j}\right)
    =\frac{\delta_{ij}}{2\left(\cc+k(r)\right)}
\end{displaymath}
as well as
\begin{displaymath}
    g_q\left(\frac{\pt}{\pt u^i},\frac{\pt}{\pt v^j}\right)=0\;.
\end{displaymath}
Moreover, the radial vector field
\begin{displaymath}
    \frac{\pt}{\pt r}=
    \frac{1}{\sqrt{\sum_{i=1}^d\left(u^i\right)^2}}
    \left(\sum_{i=1}^d u^i\frac{\pt}{\pt u^i}\right)
\end{displaymath}
satisfies at $q=(u,v)\in Q$ that
\begin{displaymath}
    g_q\left(\frac{\pt}{\pt r},\frac{\pt}{\pt r}\right)
    =\frac{1}{2\left(\cc-k(r)\right)}\;.
\end{displaymath}
In particular, the vector field $X_1$ on $Q$ defined by
\begin{displaymath}
    X_1=\sqrt{2(\cc-k(r))}\,\frac{\pt}{\pt r}
\end{displaymath}
is of unit length. As the vector field $X_1$ depends on the distance component $r$ alone, and as we can extend $X_1$ locally to an orthonormal frame $(X_1,\dots,X_{2d})$ for $TQ$,
the distance process between the two landmarks is the one-dimensional diffusion process with generator
\begin{displaymath}
    \frac{1}{2}X_1^2+\frac{1}{2}(\operatorname{div} X_1)X_1\;.
\end{displaymath}
To determine the associated It\^o stochastic differential equation, we first use
\begin{displaymath}
    \operatorname{vol}_g=
    \frac{1}{4^d\left(\cc^2-(k(r))^2\right)^{d}}\dd u^1\dots\dd u^d\dd v^1\dots \dd v^d
\end{displaymath}
to determine
\begin{align*}
    \left(\operatorname{div} X_1\right)(r)
    &=2^d\left(\cc^2-(k(r))^2\right)^{d/2}
    \frac{\pt}{\pt r}\left(
    2^{-d+1/2}\left(\cc^2-(k(r))^2\right)^{-d/2}\left(\cc-k(r)\right)^{1/2}
    \right)\\
    &=\frac{((2d-1)k(r)-\cc)k'(r)}
    {\sqrt{2\left(\cc-k(r)\right)}(\cc+k(r))}\;.
\end{align*}
We further observe that we still have, for $f\in C^2(Q)$,
\begin{displaymath}
    X_1(X_1(f))=2(\cc-k)\frac{\pt^2 f}{\pt r^2} - \frac{\pt k}{\pt r}\frac{\pt f}{\pt r}\;,
\end{displaymath}
and we compute
\begin{displaymath}
    \frac{1}{2}\left((\operatorname{div} X_1)(r)\sqrt{2(\cc-k(r))}-k'(r)\right)
    =\frac{((d-1)k(r)-\cc)k'(r)}{\cc+k(r)}\;.
\end{displaymath}
Thus, the distance process $(r_t)_{t\in[0,\zeta)}$ between the two landmarks solves the It\^o stochastic differential equation
\begin{equation}\label{eq:distanceSDE}
  \db r_t =\sigma(r_t)\dd B_t+b(r_t)\dd t\;,
\end{equation}
where $(B_t)_{t\geq 0}$ is a one-dimensional standard Brownian motion, and with the diffusivity $\sigma\colon [0,\infty)\to\R$ as well as the drift $b\colon[0,\infty)\to\R$ given by
\begin{equation}\label{eq:sigma-b}
    \sigma(r)=\sqrt{2(\cc-k(r))}
    \quad\text{and}\quad
    b(r)=\frac{((d-1)k(r)-\cc)k'(r)}{\cc+k(r)}\;.
\end{equation}
Note that this is consistent with the It\^o stochastic differential equation~\eqref{eq:distin1D} derived for $d=1$.

\subsection{Collision analysis for two landmarks}\label{sec:collision}

In the previous section, we derived the It\^o stochastic differential equation which governs the dynamics of the distance process $(r_t)_{t\in[0,\zeta)}$ between two landmarks in $\R^d$ induced by the Riemannian Brownian motion on landmark configurations consisting of exactly two landmarks. 
Studying if the two landmarks collide and aiming for a classification which depends on the choice of the kernel $K$ as well as the dimension $d\geq 1$ then amounts to analysing the one-dimensional diffusion process $(r_t)_{t\in[0,\zeta)}$ near the singularity at zero. Indeed, the landmarks collide if and only if the distance process hits zero. For our analysis, we follow the classification for singular points of one-dimensional diffusion processes by Cherny and Engelbert~\cite{SSDE}.

This classification depends on the asymptotic behaviour of the coefficients given in~(\ref{eq:sigma-b}) of the stochastic differential equation \eqref{eq:distanceSDE} near zero, and hence on the small-distance asymptotics of the kernel. Motivated by our main examples, Gaussian and Sobolev kernels, we assume that there exist real-valued constants $D,\gamma>0$ such that, as $r\downarrow 0$,
\begin{equation}\label{eq:expgamma}
    k(0)-k(r)=Dr^\gamma +o\left(r^\gamma\right)\;.
\end{equation}
To reduce the notational overhead, we write $f(r)\sim g(r)$ as $r\downarrow 0$ if there exists a non-zero constant $C\in\R$ such that, as $r\downarrow 0$,
\begin{displaymath}
    f(r)=g(r)(C+o(1))\;.
\end{displaymath}

To start off analysing the singularity at zero for the one-dimensional diffusion process $(r_t)_{t\in[0,\zeta)}$, which is the unique strong solution to the It\^o stochastic differential equation~\eqref{eq:distanceSDE}, we remark that, as $r\downarrow 0$,
\begin{displaymath}
    \frac{1+|b(r)|}{\left(\sigma(r)\right)^2}\sim r^{-\gamma}\;.
\end{displaymath}
For $0<\gamma<1$, the function above is locally integrable near zero, and according to~\cite[Definition~2.3]{SSDE}, zero is then a regular point of~\eqref{eq:distanceSDE}. 
In this case, as a consequence of~\cite[Theorem~2.11]{SSDE}, the diffusion process $(r_t)_{t\in[0,\zeta)}$ hits zero in finite time with positive probability, meaning that landmarks collide with positive probability.

The remainder of this subsection is devoted to the other case, $\gamma\geq 1$, where 
\begin{displaymath}
    \frac{1+|b|}{\sigma^2}\not\in L_{\rm loc}^1(0)\;. 
\end{displaymath}
According to~\cite[Definition~2.3]{SSDE}, zero is then a singular point of~\eqref{eq:distanceSDE}. For sufficiently small $a>0$, we have
\begin{equation}\label{eq:locintegrable}
    \frac{1+|b|}{\sigma^2}\in L_{\rm loc}^1((0,a])\;.
\end{equation}
Hence, we can proceed with the classification of singularities by Cherny and Engelbert~\cite{SSDE}, which is very well summarised on~\cite[p.~39]{SSDE}. Throughout, we fix $a>0$ such that~\eqref{eq:locintegrable} is satisfied.

In the first step of the classification process, we need to consider the function $\rho\colon (0,a]\to\R$ defined by
\begin{displaymath}
    \rho(r)=\exp\left(\int_r^a\frac{2b(y)}{\left(\sigma(y)\right)^2}\dd y\right)\;.
\end{displaymath}
We compute that
\begin{displaymath}
    \frac{2b(y)}{(\sigma(y))^2}
    =\frac{((d-1)k(y)-\cc)k'(y)}{\cc^2-\left(k(y)\right)^2}
    =\frac{(d-2)k'(y)}{2(\cc-k(y))}-\frac{d\,k'(y)}{2(\cc +k(y))}\;.
\end{displaymath}
By employing the change of variables $z=k(y)$ and subject to $a>0$ being sufficiently small, we further obtain
\begin{align*}
    \int_r^a\frac{2b(y)}{(\sigma(y))^2}\dd y
    &=\int_r^a\frac{(d-2)k'(y)}{2(\cc-k(y))}\dd y-\int_r^a\frac{d\,k'(y)}{2(\cc +k(y))}\dd y\\
    &=\int_{k(r)}^{k(a)}\frac{d-2}{2(\cc-z)}\dd z-\int_{k(r)}^{k(a)}\frac{d}{2(\cc + z)}\dd z\\
    &=\frac{1}{2}\left[(2-d)\log(\cc-z)-d\log(\cc+z)\right]_{k(r)}^{k(a)}\\
    &=\frac{1}{2}\left((2-d)\log\left(\frac{\cc-k(a)}{\cc-k(r)}\right)-d\log\left(\frac{\cc+k(a)}{\cc+k(r)}\right)\right)\;,
\end{align*}
which yields
\begin{displaymath}
    \rho(r)=
    \left(\frac{\cc-k(a)}{\cc-k(r)}\right)^{1-d/2}\left(\frac{\cc+k(a)}{\cc+k(r)}\right)^{-d/2}\;.
\end{displaymath}
From~\eqref{eq:expgamma} it then follows directly that, as $r\downarrow 0$,
\begin{equation}\label{eq:rhoasymp}
    \rho(r)\sim r^{(d/2-1)\gamma}\;,
\end{equation}
which shows that
\begin{displaymath}
    \int_0^a\rho(r)\dd r
    \begin{cases}
        =\infty & \text{if }d=1\text{ and }\gamma\geq 2\;,\\
        <\infty & \text{otherwise}\;.
    \end{cases}
\end{displaymath}
We further deduce that, as $r\downarrow 0$,
\begin{displaymath}
    \frac{1+|b(r)|}{\rho(r)(\sigma(r))^2}\sim
    r^{-(d/2-1)\gamma}r^{-\gamma}\sim r^{-d\gamma/2}\;,    
\end{displaymath}
implying that
\begin{displaymath}
    \int_0^a\frac{1+|b(r)|}{\rho(r)(\sigma(r))^2}\dd r
    \begin{cases}
        <\infty & \text{if }d=1\text{ and }1\leq\gamma<2\;,\\
        =\infty & \text{otherwise}\;.
    \end{cases}
\end{displaymath}

To complete the classification for $d=1$ and $1\leq\gamma<2$, we observe that in this case, as $r\downarrow 0$,
\begin{displaymath}
    \frac{|b(r)|}{(\sigma(r))^2}
    =\frac{|\cc k'(r)|}{2(\cc^2-(k(r))^2)}\sim r^{-1}\;,
\end{displaymath}
from which we conclude
\begin{displaymath}
    \int_0^a \frac{|b(r)|}{(\sigma(r))^2}\dd r=\infty\;.
\end{displaymath}
Thus, according to the result~\cite[Theorem~2.12]{SSDE}, for $d=1$ and $1\leq\gamma<2$, the singularity of~\eqref{eq:distanceSDE} at zero is of type 2. In particular, the associated distance process $(r_t)_{t\in[0,\zeta)}$ hits zero with positive probability, that is, the two landmarks collide with positive probability in finite time.

The remaining classification steps use the function $s\colon (0,a]\to\R$ defined by
\begin{displaymath}
    s(r)=
    \begin{cases}
        \int_0^r\rho(y)\dd y & \text{if } \int_0^a\rho(y)\dd y<\infty\;,\\
        \int_a^r\rho(y)\dd y & \text{if } \int_0^a\rho(y)\dd y=\infty\;.
    \end{cases}
\end{displaymath}
From~\eqref{eq:rhoasymp} we obtain that, as $r\downarrow 0$,
\begin{displaymath}
    s(r)\sim
    \begin{cases}
        \log(r) & \text{if }d=1\text{ and }\gamma=2\;,\\
        r^{(d/2-1)\gamma+1} & \text{otherwise}\;.        
    \end{cases}
\end{displaymath}
Therefore, for $d=1$ and $\gamma=2$, we have, as $r\downarrow 0$,
\begin{displaymath}
    \frac{1+|b(r)|}{\rho(r)(\sigma(r))^2}s(r)
    \sim r^{-1}\log(r)\;,
\end{displaymath}
whilst in all other cases, we have, as $r\downarrow 0$,
\begin{displaymath}
    \frac{1+|b(r)|}{\rho(r)(\sigma(r))^2}s(r)
    \sim r^{-d\gamma/2}r^{(d/2-1)\gamma+1}
    \sim r^{1-\gamma}\;.
\end{displaymath}
It follows that, irrespective of the dimension $d\geq 1$,
\begin{displaymath}
    \int_0^a \frac{1+|b(r)|}{\rho(r)(\sigma(r))^2}s(r) \dd r
    \begin{cases}
        <\infty & \text{if } 1\leq\gamma<2\;,\\
        =\infty & \text{if } \gamma\geq 2\;.
    \end{cases}
\end{displaymath}
At this stage, we deduce that if $d\geq 2$ and $1\leq\gamma<2$ then by~\cite[Theorem~2.13]{SSDE} the singularity of~\eqref{eq:distanceSDE} at zero is of type 1, which particularly implies that in this case the two landmarks collide with positive probability in finite time. For $d=1$ and $\gamma\geq 2$, it is a consequence of~\cite[Theorem~2.17]{SSDE} that the singularity of~\eqref{eq:distanceSDE} at zero is then of type 5. Hence, in this case any solution to~\eqref{eq:distanceSDE} started at a non-zero distance is strictly positive. Together with the argument presented in the following section, this implies that the associated landmark Brownian motion exists for all times.

To conclude the classification of the singularity at zero for the distance process $(r_t)_{t\in[0,\zeta)}$, it remains to observe that, for $d\geq 2$ and $\gamma\geq 2$, as $r\downarrow 0$,
\begin{displaymath}
    \frac{s(r)}{\rho(r)(\sigma(r))^2}
    \sim r^{(d/2-1)\gamma+1}r^{-(d/2-1)\gamma}r^{-\gamma}
    \sim r^{1-\gamma}\;,
\end{displaymath}
which yields that, for $d\geq 2$ and $\gamma\geq 2$,
\begin{displaymath}
    \int_0^a \frac{s(r)}{\rho(r)(\sigma(r))^2}\dd r = \infty\;.
\end{displaymath}
Hence, by~\cite[Theorem~2.15]{SSDE}, in the case $d\geq 2$ and $\gamma\geq 2$, the singularity of~\eqref{eq:distanceSDE} at zero is of type 4. Whilst this together with Section~\ref{sec:escape} still establishes long-time existence of the associated landmark Brownian motion, the type of singularity detected implies that with positive probability the two landmarks draw arbitrarily close together, if measured with respect to the Euclidean metric.

A summary of the results derived is provided in Section~\ref{sec:application} below.

\subsection{Ruling out escape to infinity}\label{sec:escape}

Besides collision of landmarks, the only other source of Brownian incompleteness is escape to infinity, which we investigate next. 
We show that escape to infinity cannot occur before collision and therefore plays no role in the analysis of Brownian completeness or incompleteness in our setting.
The presented proof requires that the kernel $k$ and its derivative $k'$ are Lipschitz continuous away from zero and that $k$ vanishes at infinity, as shall be assumed herewith. 
As before, $(q_t)_{t\in[0,\zeta)}$ denotes Brownian motion on the space of two landmarks in $\mathbb R^d$. 
Moreover, we write $r(q)$ for the inter-landmark distance of $q\in Q$.
Then, the escape-to-infinity time $S$ and the collision time $T$ are defined as
\begin{align*}
S_N &= \inf \{t>0:\|q_t\|_{\R^{nd}}\geq N\}\;,
&
S &= \lim_{N\to \infty} S_N\;,
\\
T_N &= \inf \{t>0:r(q_t)\leq 1/N\}\;,
&
T &= \lim_{N\to \infty} T_N\;.
\end{align*}
We claim that $S\geq T$, which we
establish by proving the equivalent result that $S\geq T_N$ for all $N \in \mathbb N$. 
To this end, we truncate the stochastic differential equation \eqref{eq:RiemBM} for $(q_t)_{t\in[0,\zeta)}$ using a Lipschitz function $\varphi_N\colon\mathbb R\to[0,1]$ such that $\varphi_N(r)=0$ for $r\leq1/(N+1)$ and $\varphi_N(r)=1$ for $r\geq 1/N$ to obtain 
\begin{align}\label{eq:truncSDE}
\begin{aligned}
    \dd q^{(N),i}_t &= \varphi_N\left(r\left(q^{(N)}_t\right)\right)\sqrt{K\left(q^{(N)}_t\right)}^i
    \dd W_t \\
    &\qquad - \frac{1}{2}\varphi_N\left(r\left(q^{(N)}_t\right)\right)\sum_{\ell,m=1}^n
    K\left(q^{(N)}_t\right)^{\ell m}\Gamma\left(q^{(N)}_t\right)^i_{\ell m}\dd t\;. 
\end{aligned}
\end{align}
This truncated stochastic differential equation is well-posed because its coefficients are Lipschitz continuous, as we show below. 
Hence, $(q^{(N)}_t)_{t\in[0,\zeta_N)}$ does not escape to infinity in finite time. 
Moreover, $(q^{(N)}_t)_{t\in[0,\zeta_N)}$ coincides with the Brownian motion $(q_t)_{t\in[0,\zeta)}$ on the stochastic interval $[0,T_N]$. 
Consequently, $(q_t)_{t\in[0,\zeta)}$ does not escape to infinity before $T_N$, that is, $S\geq T_N$ for all $N\in\N$, as claimed. 

We now show that the diffusivity coefficient in~(\ref{eq:truncSDE}) is Lipschitz continuous on the set $Q_N$ of all landmark configurations $(x,y)$ with inter-landmark distance $r=\|x-y\|_{\R^d}\geq 1/N$. As the scalar kernel $k$ is assumed to be Lipschitz on $[1/N,\infty)$, the cometric $g^{-1}$ is Lipschitz continuous on $Q_N$. Moreover, as the scalar kernel $k$ by assumption extends continuously to the compact set $[1/N,\infty]$, the set $g^{-1}(Q_N)$ has compact closure, namely, the set $g^{-1}(Q_N)$ itself together with the matrix
\begin{displaymath}
    \begin{pmatrix}
        k(0)&0\\0&k(0)
    \end{pmatrix}
    \otimes I_d\;.
\end{displaymath}
Thus, all of these matrices are positive definite. Taking the square root of a symmetric positive definite matrix is smooth by the implicit function theorem or, more generally, because the functional calculus is real analytic, see~\cite{bauer2022smooth}. In particular, the matrix square root is Lipschitz continuous on compacts. Consequently, $g^{-1/2}$ is Lipschitz continuous on $Q_N$. 

We next show that the drift coefficient of~(\ref{eq:truncSDE}) is Lipschitz continuous on $Q_N$. The derivative $k'$ of the scalar kernel is assumed to be Lipschitz continuous on $[1/N,\infty)$. Therefore, the cometric $g^{-1}$ has Lipschitz continuous coordinate derivatives on $Q_N$. Matrix inversion is real analytic and hence Lipschitz continuous on compacts. Thus, the metric $g$ is Lipschitz continuous on $Q_N$. The Christoffel symbol can be written as a contraction of the metric $g$ with coordinate derivatives of the cometric $g^{-1}$. Therefore, the Christoffel symbol is Lipschitz continuous on $Q_N$. Taken together, this implies the Lipschitz continuity of the drift in~(\ref{eq:truncSDE}) on $Q_N$.

This concludes the proof that $S\geq T$, that is, escape to infinity cannot occur before collision. An important consequence is that Brownian completeness follows as soon as collisions are ruled out.

\subsection{Summary}\label{sec:application}

The long-term behaviour of Riemannian Brownian motion on the configuration space of two landmarks in $\mathbb R^d$ depends on whether the ambient space has dimension $d=1$ or $d\geq 2$, and on whether the near-zero asymptotics \eqref{eq:expgamma} of the kernel are given by $\gamma<2$ or $\gamma\geq 2$. Our characterisation follows \cite{SSDE} and is well described in terms of the two hitting times $T_a=\inf\{t\geq 0:r_t=a\}$ and $T_{0,a}=\min\{T_0,T_a\}$. 

\begin{itemize}
\item If $\gamma<2$, then the two landmarks collide with positive probability. More specifically, for $\gamma \in (0,1)$, the origin is a regular point of the distance process $(r_t)_{t\in[0,\zeta)}$, which implies that $\E[T_{0,a}]<\infty$ and $\mathbb{P}(X_{T_{0,a}}=0)>0$. For $\gamma \in [1,2)$, it is a singular point of type 2 if $d=1$ and of type 1 if $d\geq 2$. This means that, subject to $r_0\in[0,a]$, in the case $d=1$, there exists a unique solution to~\eqref{eq:distanceSDE} up to $T_a$ such that $\E[T_a]<\infty$ and $\mathbb{P}(\text{there exists }t\leq T_a\text{ such that }X_t=0)>0$, whereas if $d\geq 2$, there exists a unique solution to~\eqref{eq:distanceSDE} which is defined up to $T_{0,a}$ and has $\E[T_{0,a}]<\infty$ as well as $\mathbb{P}(X_{T_{0,a}}=0)>0$.

\item If $\gamma\geq 2$ and $d=1$, then the landmark Brownian motion exists for all times.
The singular point at zero is of type 5, meaning that for $r_0>0$, any solution to~\eqref{eq:distanceSDE} is strictly positive, and subject to $r_0\in(0,a]$ there exists a unique solution defined up to $T_a$ and $T_a<\infty$ $\mathbb{P}$-a.s. 

\item If $\gamma\geq 2$ and $d\geq 2$, then the landmark Brownian motion exists for all times. The singularity at zero is of type 4, which implies that as long as $r_0>0$ any solution to~\eqref{eq:distanceSDE} is strictly positive, and for $r_0\in(0,a]$ there exists a unique solution defined up to $T_a$ where $\mathbb{P}(T_a=\infty)>0$ as well as $\lim_{t\to \infty} r_t=0$ $\mathbb{P}$-a.s. on $\{T_a=\infty\}$. In particular, the two landmarks almost surely do not collide, but with positive probability their Euclidean distance becomes arbitrarily small as time tends to infinity.
\end{itemize}

We next discuss implications for the Bessel potentials defined by~\eqref{eq:bessel}, which are also known as Sobolev kernels. 
The Bessel potential of order $\alpha>d$ in $d$ dimensions has asymptotics \eqref{eq:expgamma} with $\gamma=\min\{\alpha-d,2\}$, as can be seen from \eqref{eq:bessel-asymptotics}. 
A minor modification is needed in the case $\gamma=2$ to accommodate the logarithmic term in the asymptotics \eqref{eq:bessel-asymptotics}, but careful inspection of the arguments in Section~\ref{sec:collision} shows that the conclusion remains the same as for $\gamma=2$ without the logarithmic term. 

To summarise, the Sobolev metric of order $\alpha>d$ gives rise to a Brownian complete Riemannian manifold if $\alpha\geq d+2$ and to a Riemannian manifold which is Brownian incomplete otherwise.
Note that this is also the threshold for the reproducing kernel Hilbert space $\mathfrak X_{H_{\alpha/2}}(\R^d)$ to embed into $\mathfrak X_{C^1_b}(\R^d)$, that is, for admissibility of this space of vector fields in the terminology of \cite{younes}. 
Interestingly, we have Brownian completeness not only above this threshold but also in the critical case $\alpha=d+2$.
For the Gaussian kernel, one has both Brownian completeness and admissibility, in line with the interpretation of the Gaussian kernel as a Sobolev kernel of infinite order. 

\section{Numerical experiments}\label{sec:numerics}

We simulate Brownian motion of $n=2$ landmarks in $\R^d$ with $d\in\{1,2\}$. Subsequently, we repeat the experiment for $n>2$ landmarks to give numerical hints on the possibility of collision in cases which are not covered by our theoretical results. 
For kernels, we choose Bessel potentials \eqref{eq:bessel} with parameters $\nu:=(\alpha-d)/2 \in \{1/2, 3/2\}$ and the Gaussian kernel, as these admit explicit formulae. Specifically, the kernels we use in our simulations are given by
\begin{equation*}
k_{1/2}(r) = \e^{-r}\;, 
\qquad
k_{3/2}(r) = 2(1+r)\e^{-r}\;, 
\qquad
k_G(r) = \e^{-r^2}\;.
\end{equation*}
These kernels coincide with the ones in Section~\ref{sec:kernels} up to positive constants, which do not affect Brownian completeness or incompleteness.

According to our theoretical analysis, the kernel $k_{1/2}$ is Brownian incomplete, whereas $k_{3/2}$ and higher-order kernels including $k_G$ are Brownian complete, for $n=2$ landmarks in arbitrary dimension $d$.

For our experiments, we draw 20 sample paths, simulated from $t=0$ to $t=1$ with $10^4$ steps. The simulations are stopped if the inter-landmark distance gets very small or decreases very rapidly; this is taken as an indication for collision. Due to the time-discretisation and potential numerical instability for nearby landmarks, the detection of collisions is only indicative. However, as shown below, the numerical experiments are all in line with the theoretical predictions.

\subsection{Two landmarks}

\begin{figure}
    \centering
    \includegraphics[width=0.32\textwidth]{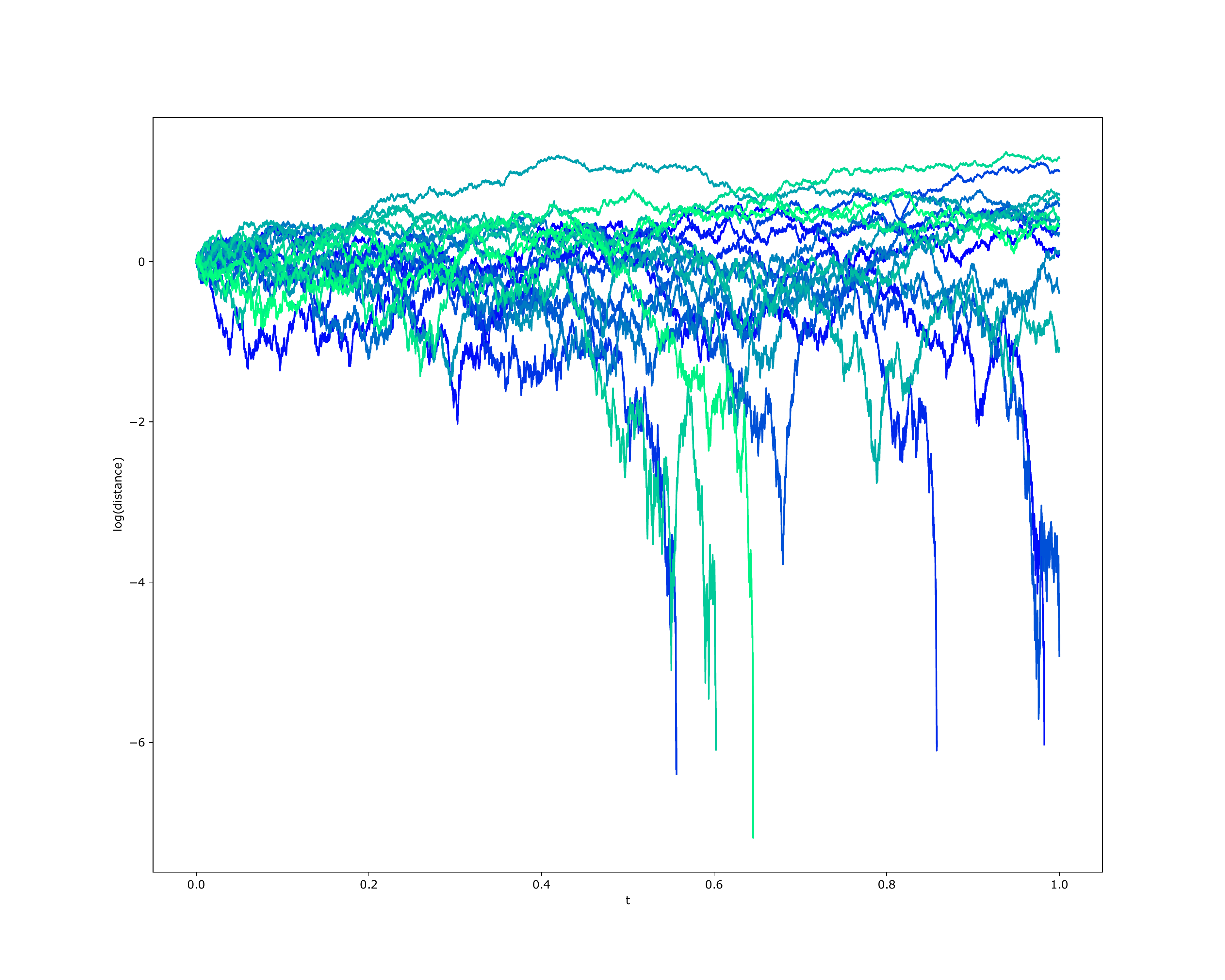}
    \includegraphics[width=0.32\textwidth]{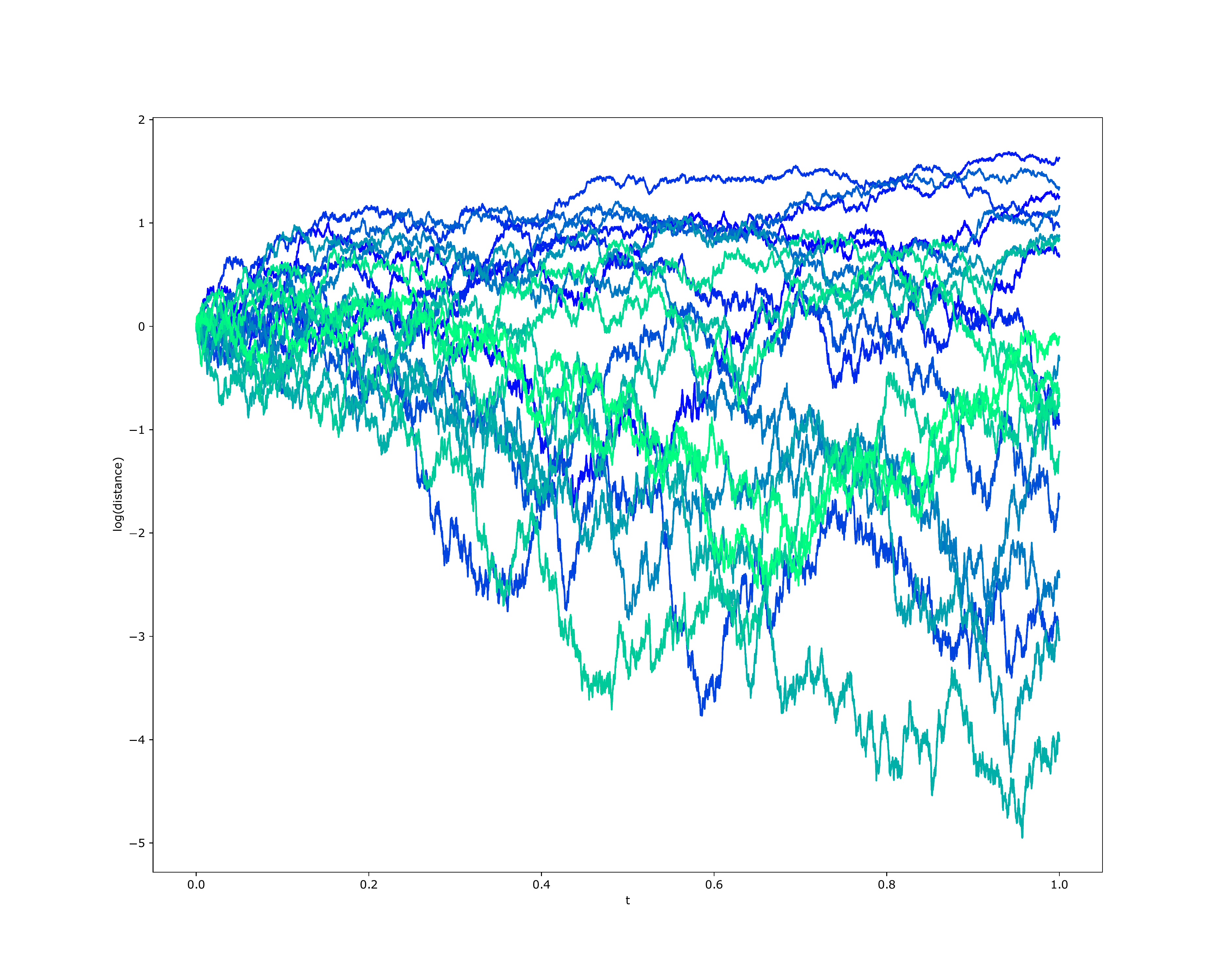}
    \includegraphics[width=0.32\textwidth]{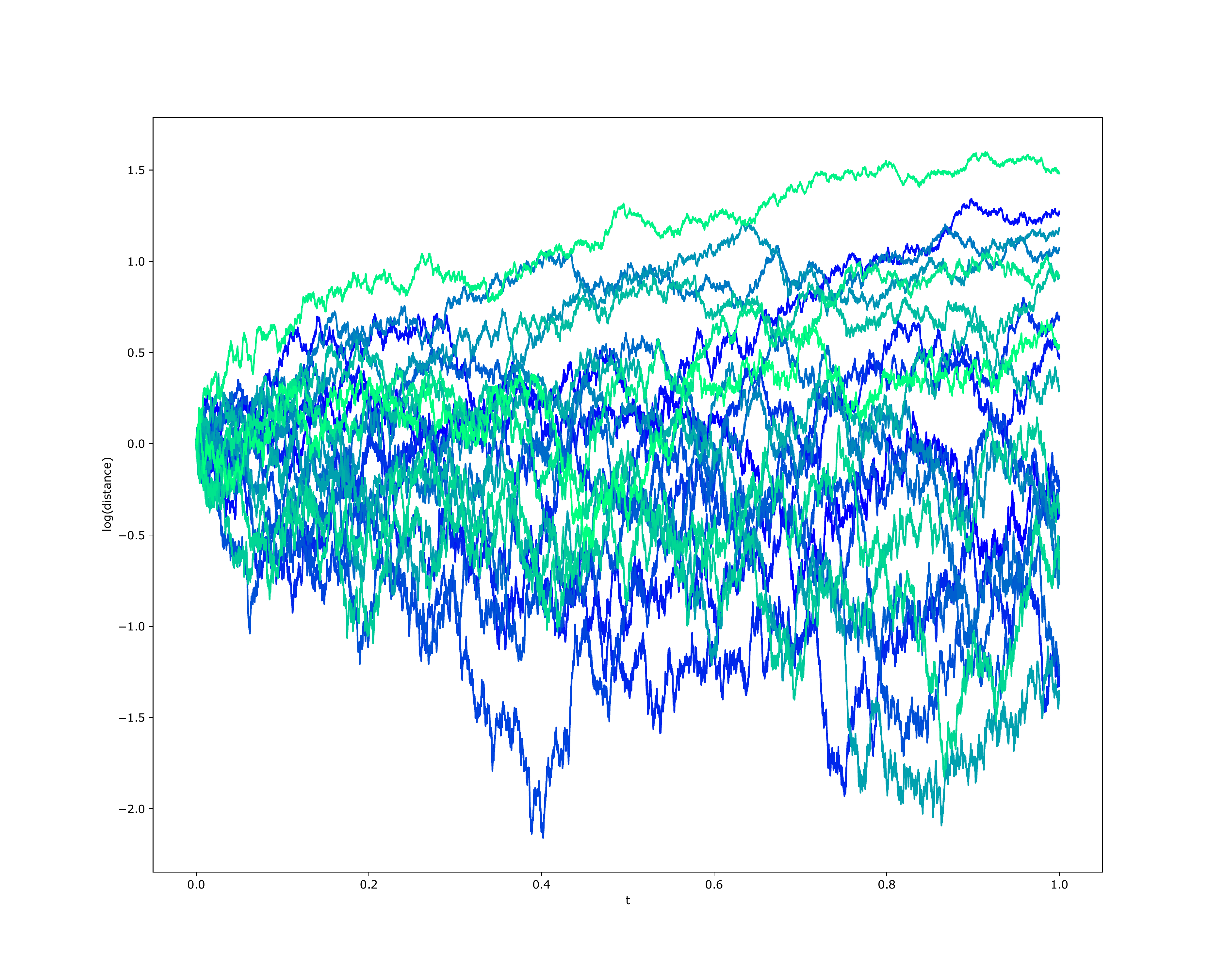}

    \includegraphics[width=0.32\textwidth]{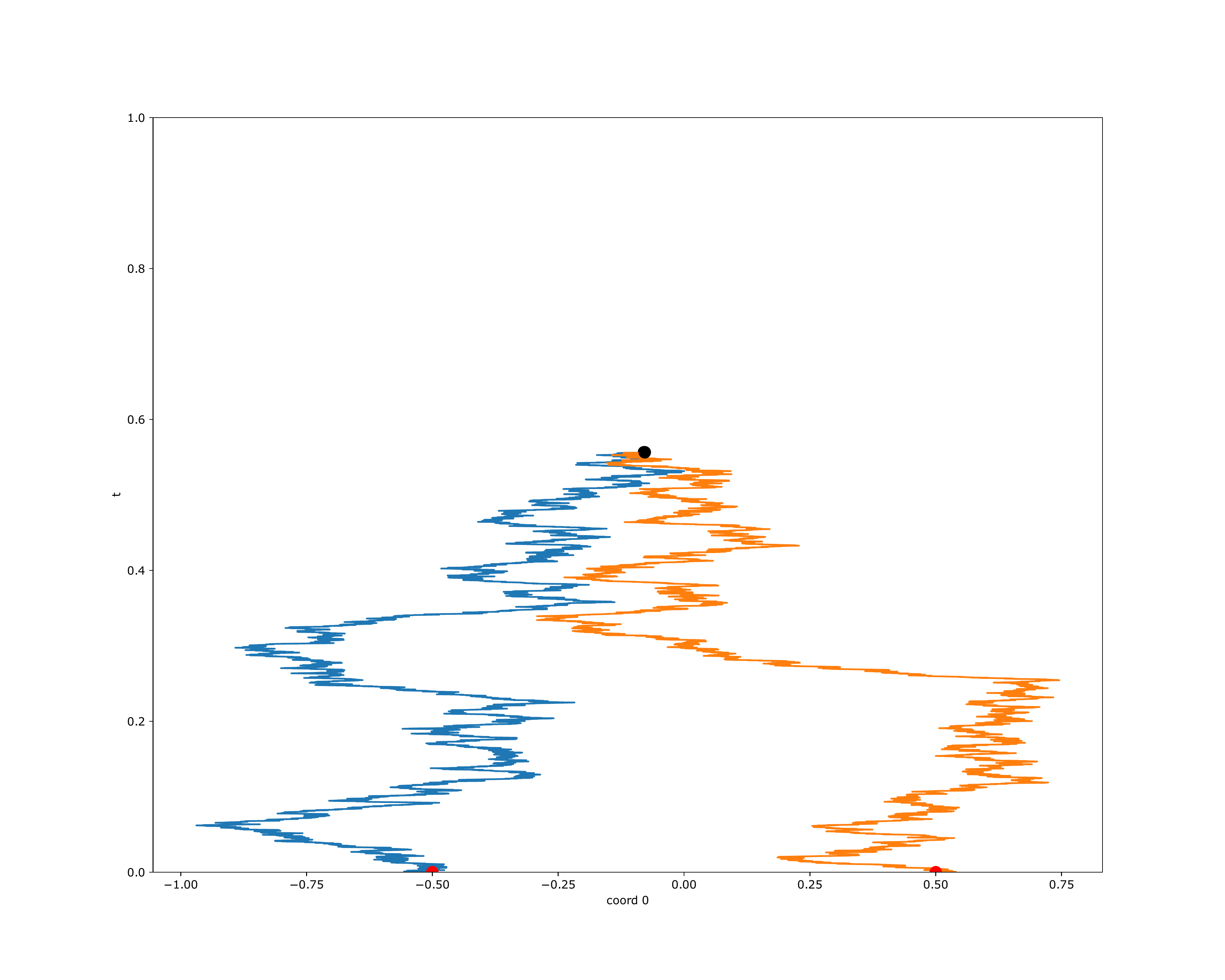}
    \includegraphics[width=0.32\textwidth]{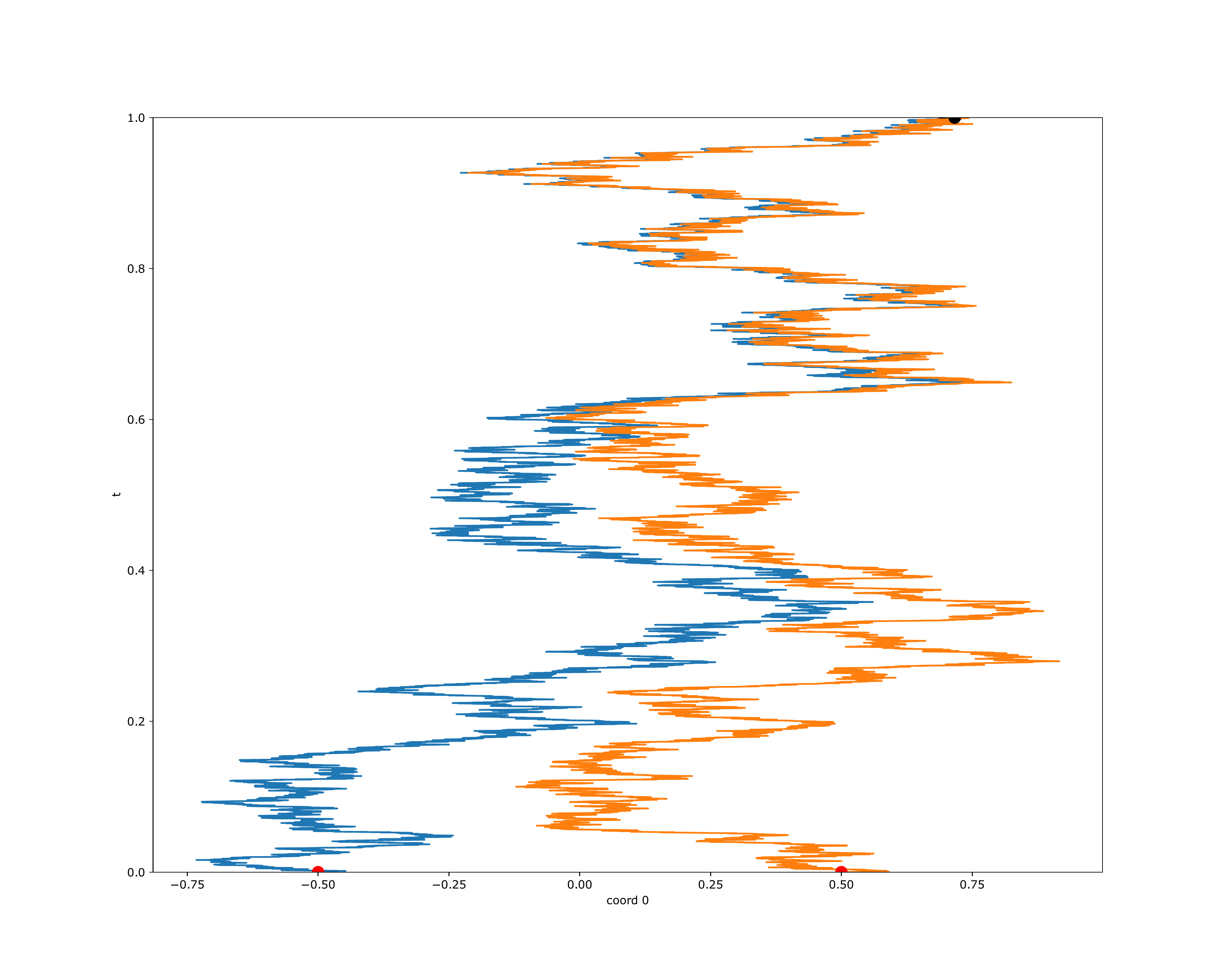}
    \includegraphics[width=0.32\textwidth]{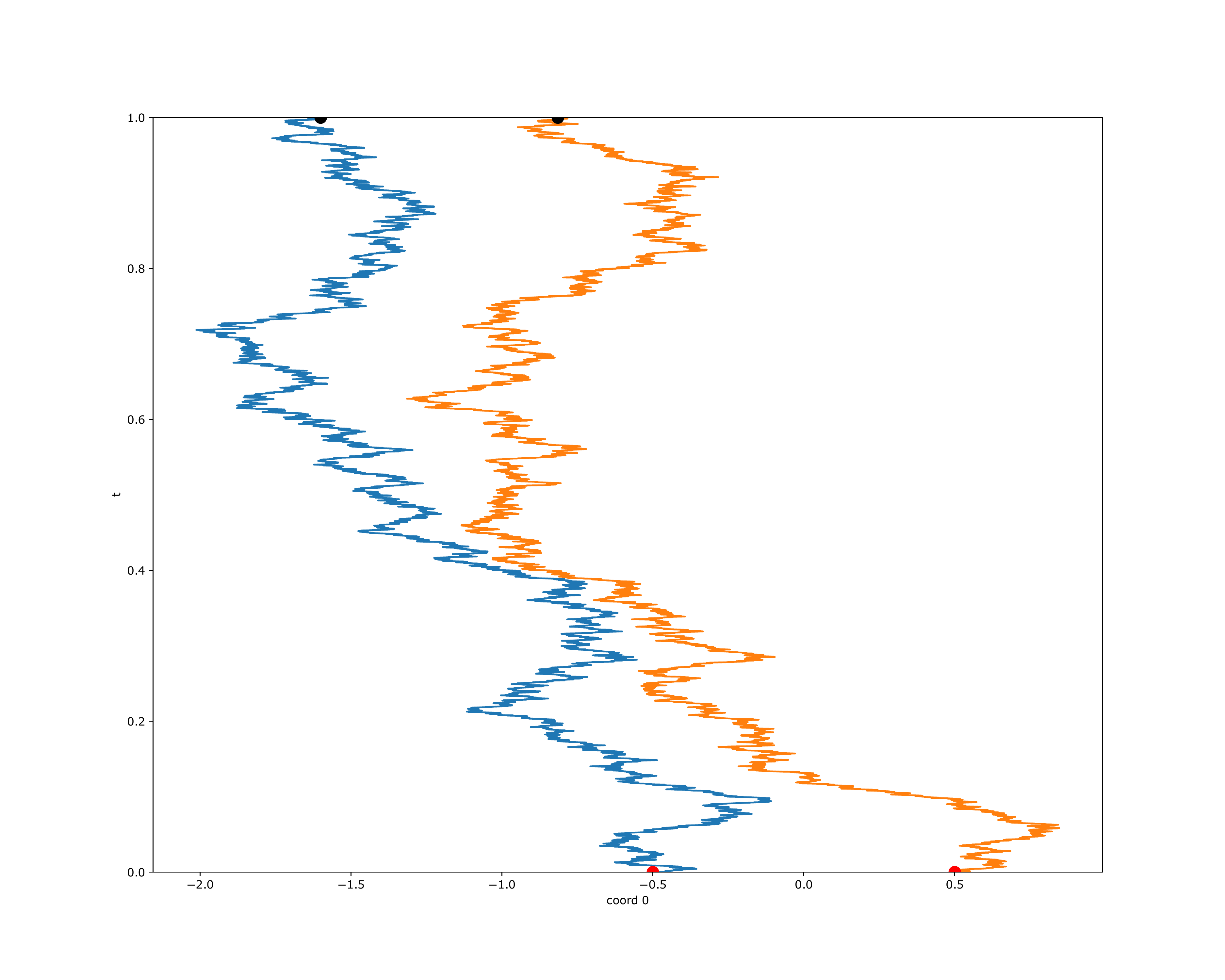}

    \caption{Results for $d=1$ and $n=2$ with 20 sample paths. 
    The first column shows the results for the the kernel $k_{1/2}$, the second column for the kernel $k_{3/2}$, and the third column for the Gaussian kernel.
    The first row shows the log-distances between the two landmarks for all sample paths, stopped if collision occurs. The second row shows the position of the two landmarks on the horizontal axis and time on the vertical axis for the sample path attaining the smallest inter-landmark distance, again stopped if the landmarks collide. 
    Whilst the landmarks temporarily come close for all kernels, a rapid decrease in the distance is observed only for the kernel $k_{1/2}$, which indicates collision.}
    \label{fig:1D}
\end{figure}

Figure~\ref{fig:1D} shows the results in dimension $d=1$ for the kernels $k_{1/2}$, $k_{3/2}$ and $k_G$. For each experiment, we sample 20 paths from the Riemannian Brownian motion and compute the distance between the two landmarks. We observe collision for $k_{1/2}$ but not for $k_{3/2}$ or $k_G$, consistent with our theoretical results. 
Higher-order Sobolev kernels lead to similar results like $k_{3/2}$, as predicted by the theory, and are not shown here.
\begin{figure}
    \centering
    \includegraphics[width=0.32\textwidth]{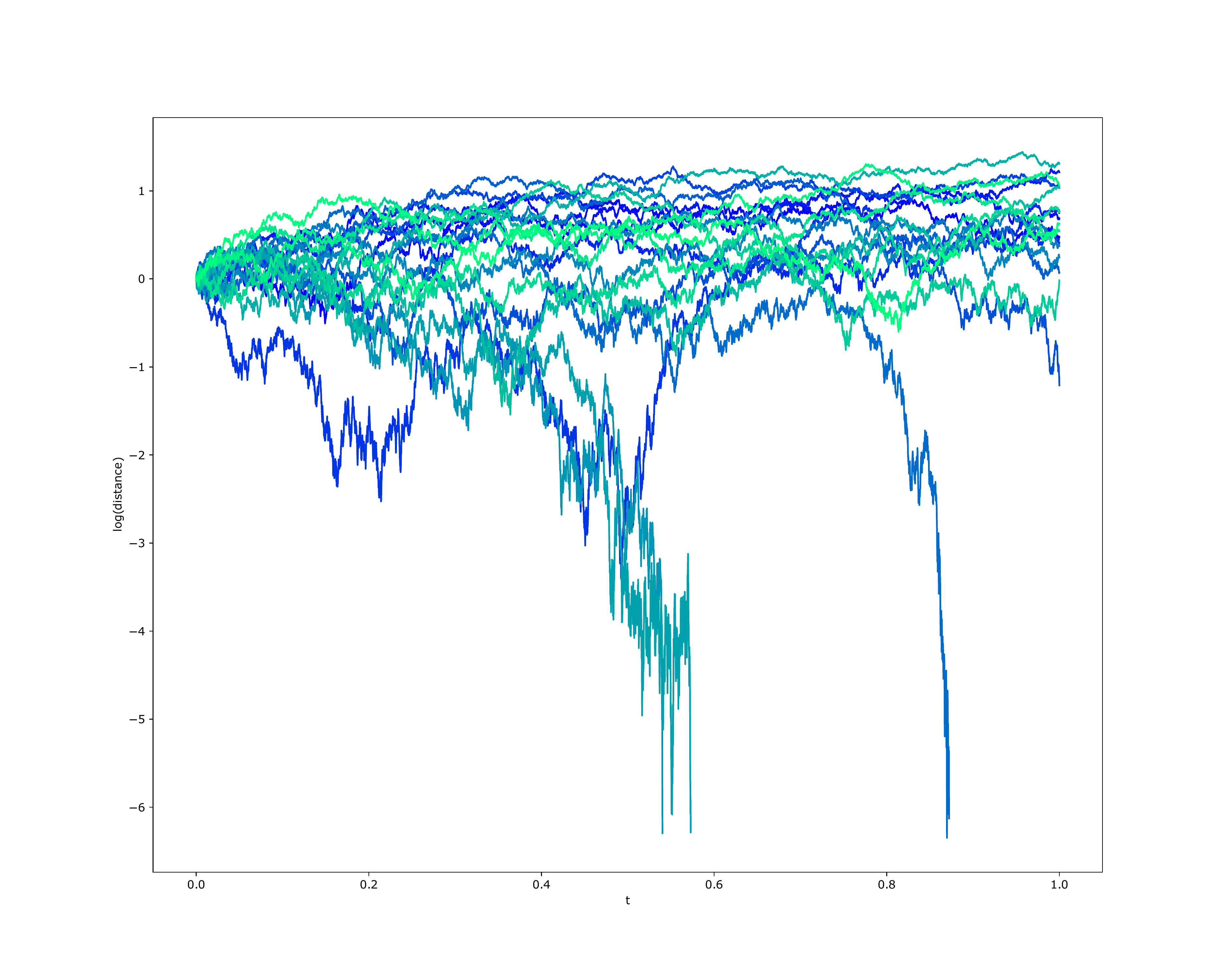}
    \includegraphics[width=0.32\textwidth]{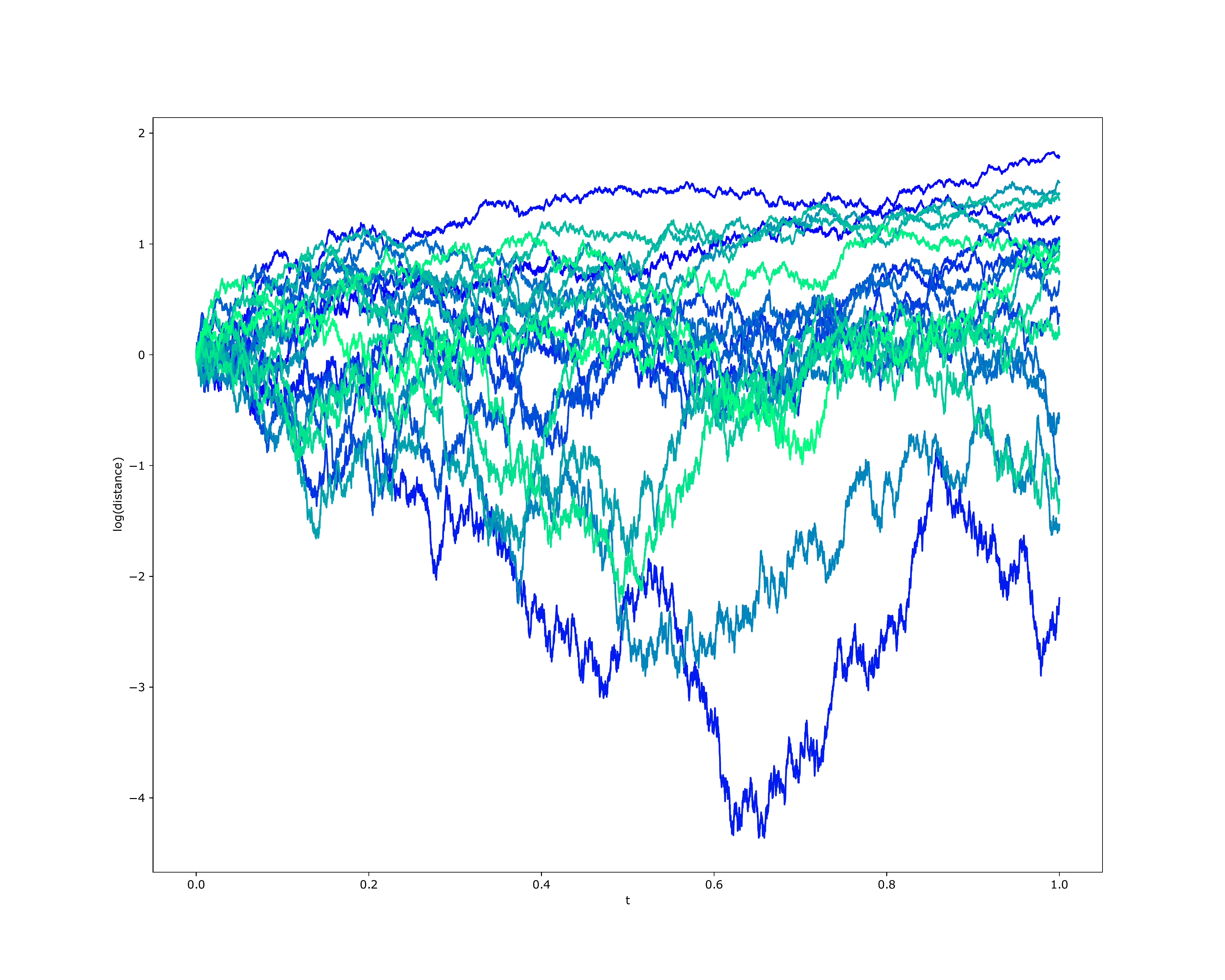}
    \includegraphics[width=0.32\textwidth]{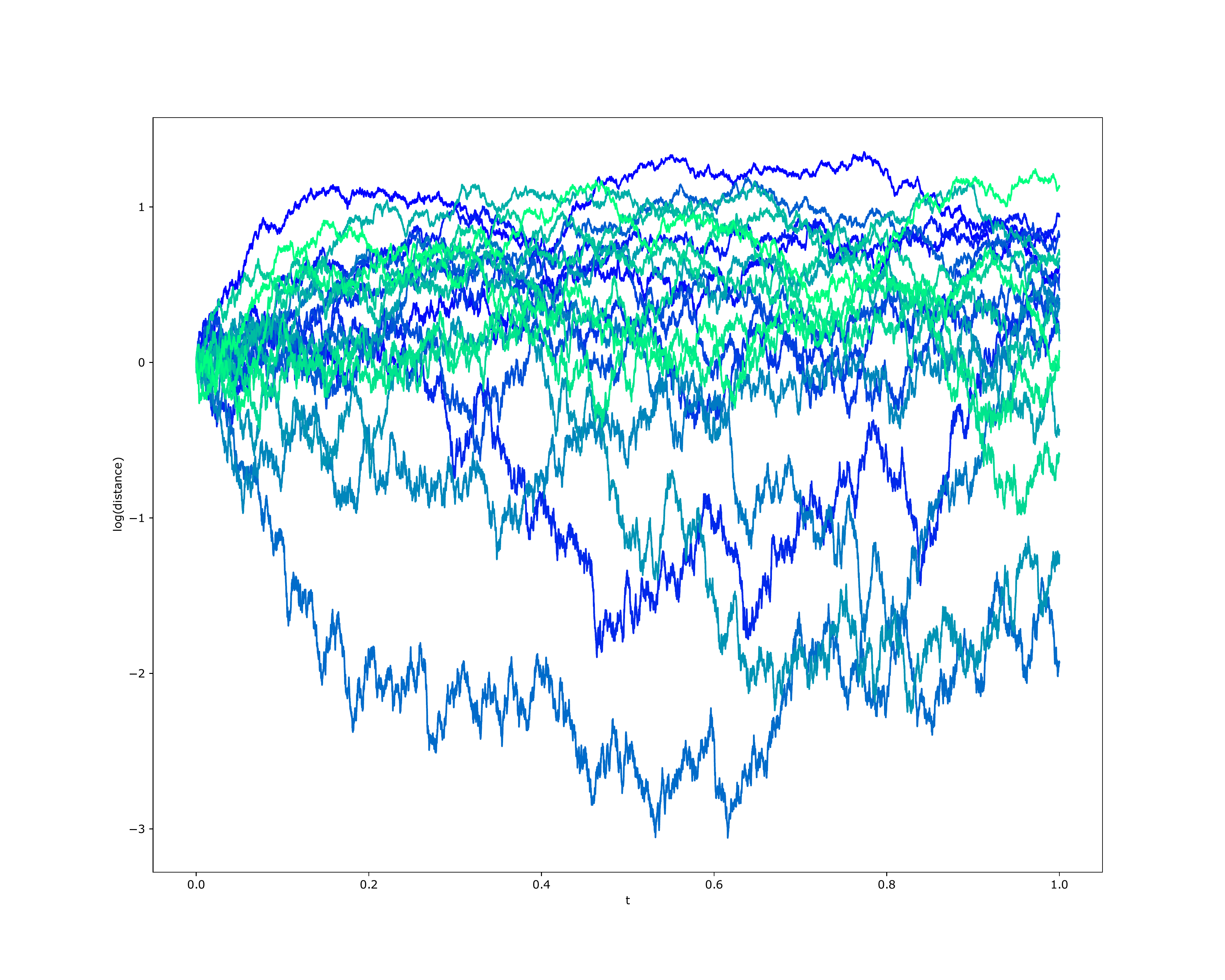}

    \includegraphics[width=0.32\textwidth]{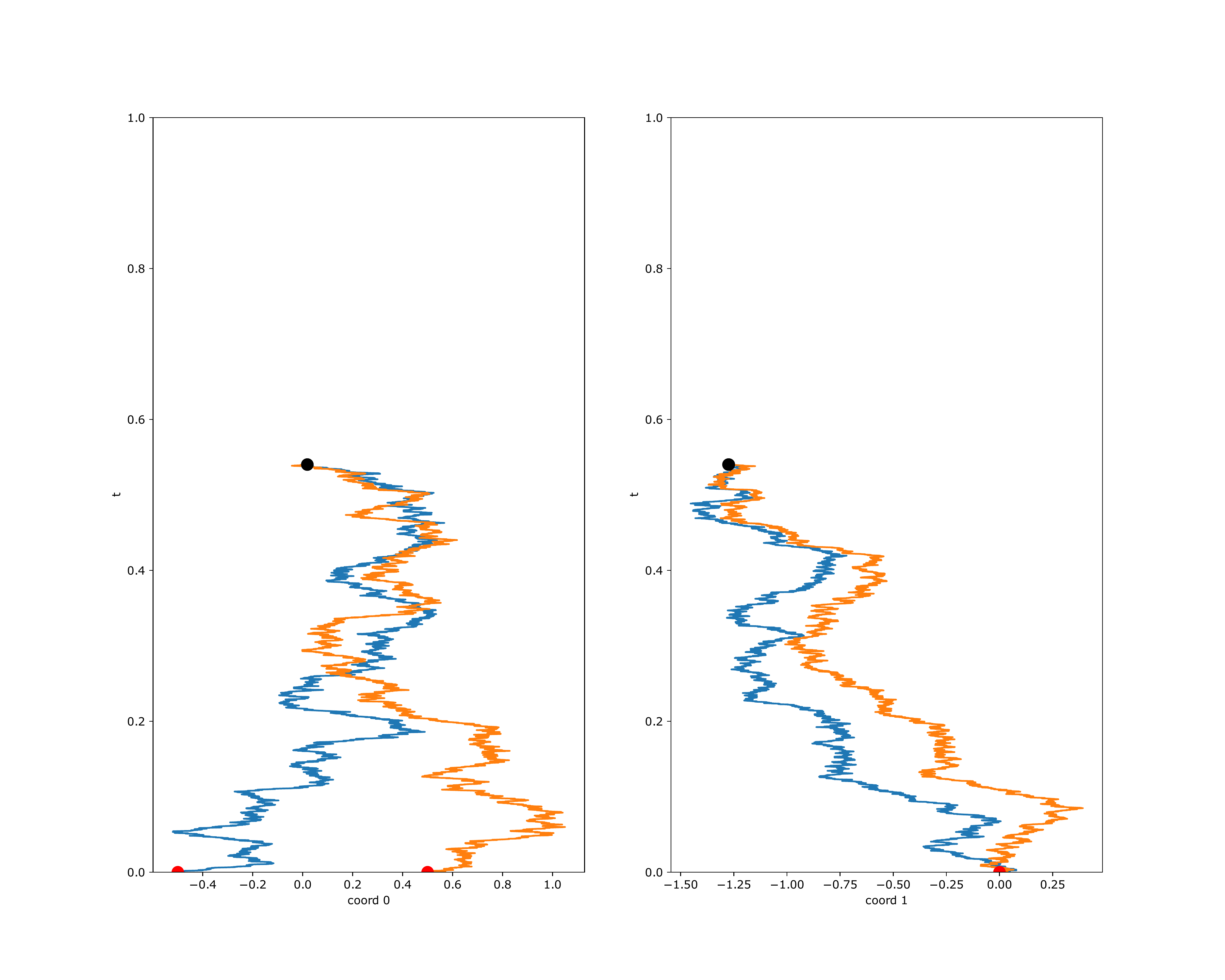}
    \includegraphics[width=0.32\textwidth]{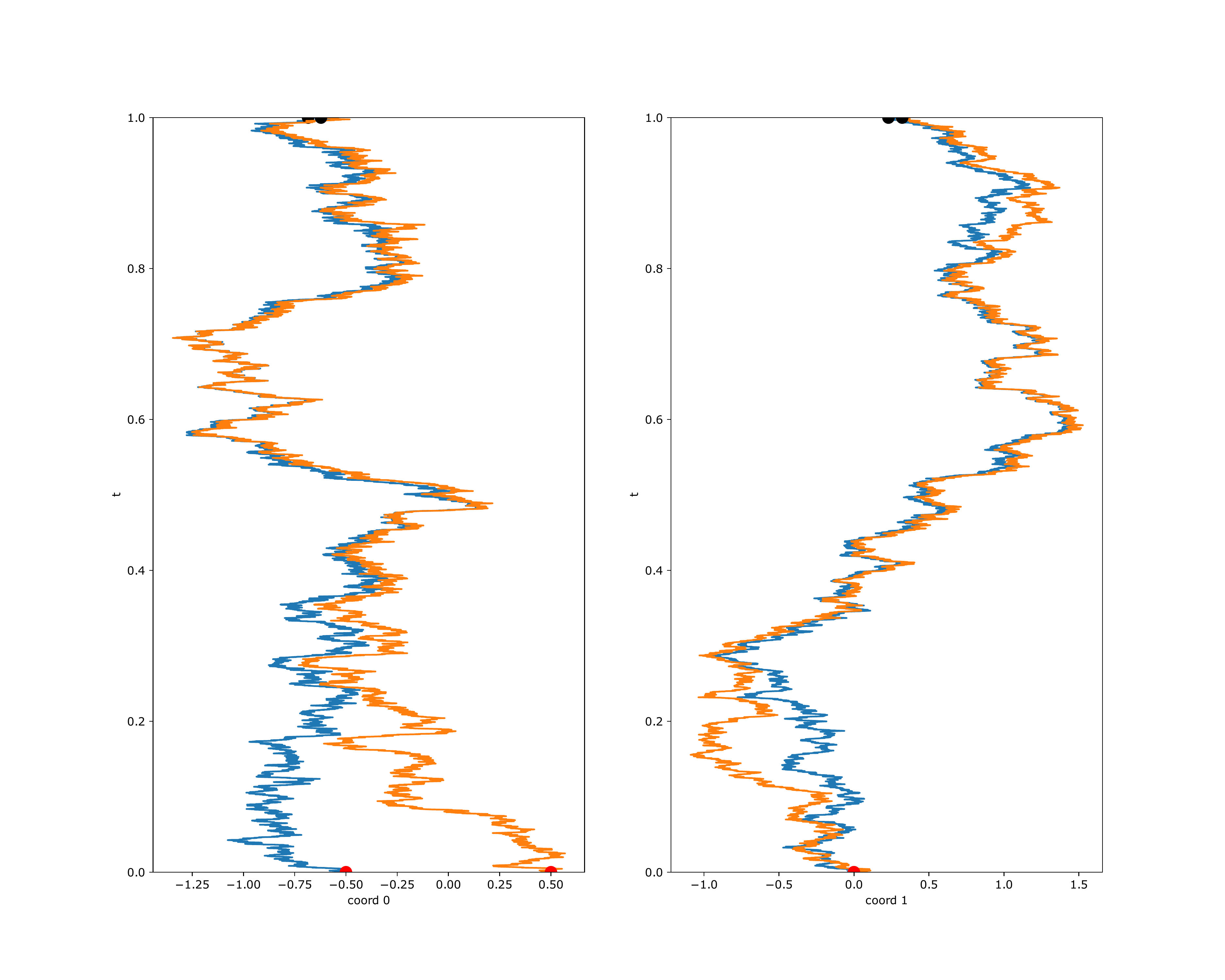}
    \includegraphics[width=0.32\textwidth]{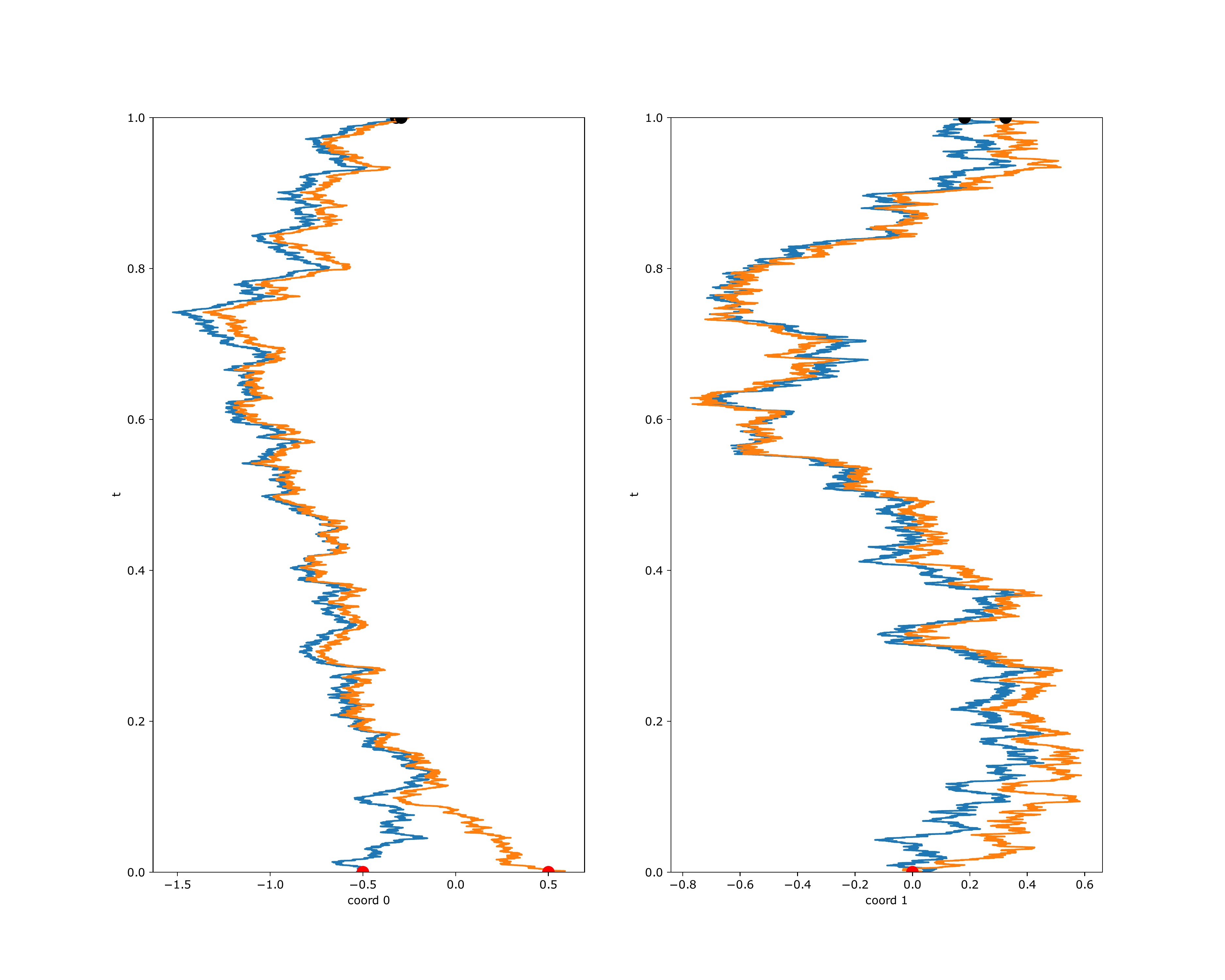}

    \caption{Results for $d=2$ and $n=2$. Setup as in Figure~\ref{fig:1D}.}
    \label{fig:2D}
\end{figure}

Figure~\ref{fig:2D} shows corresponding results in dimension $d=2$. The results are similar to dimension $d=1$, with collision observed for the kernel $k_{1/2}$ only. Note that all coordinates of the landmarks must coincide simultaneously for a collision to occur. The two coordinates are displayed separately in the plots in the second row. 

\subsection{More than two landmarks}
\begin{figure}
    \centering
    \includegraphics[width=0.32\textwidth]{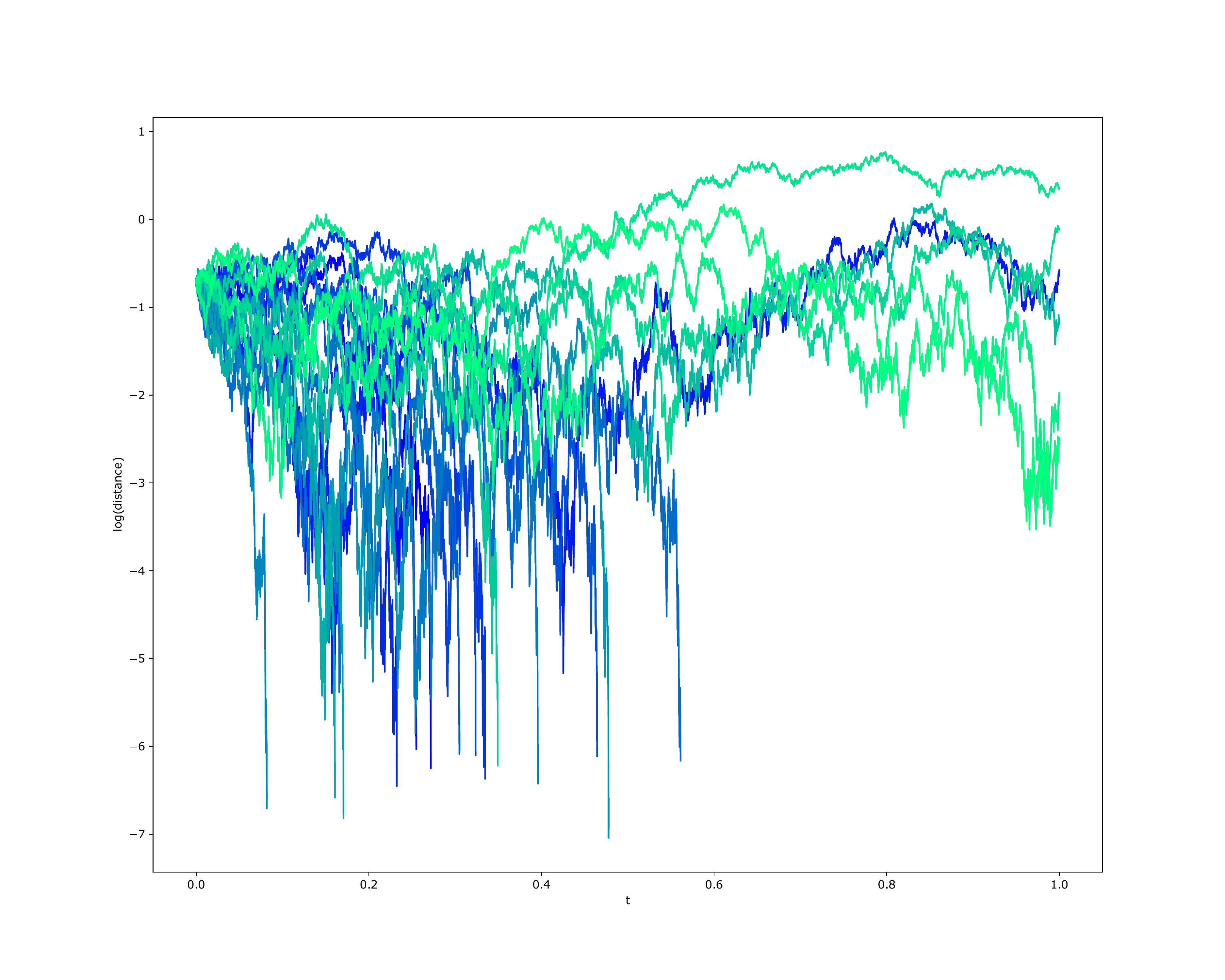}
    \includegraphics[width=0.32\textwidth]{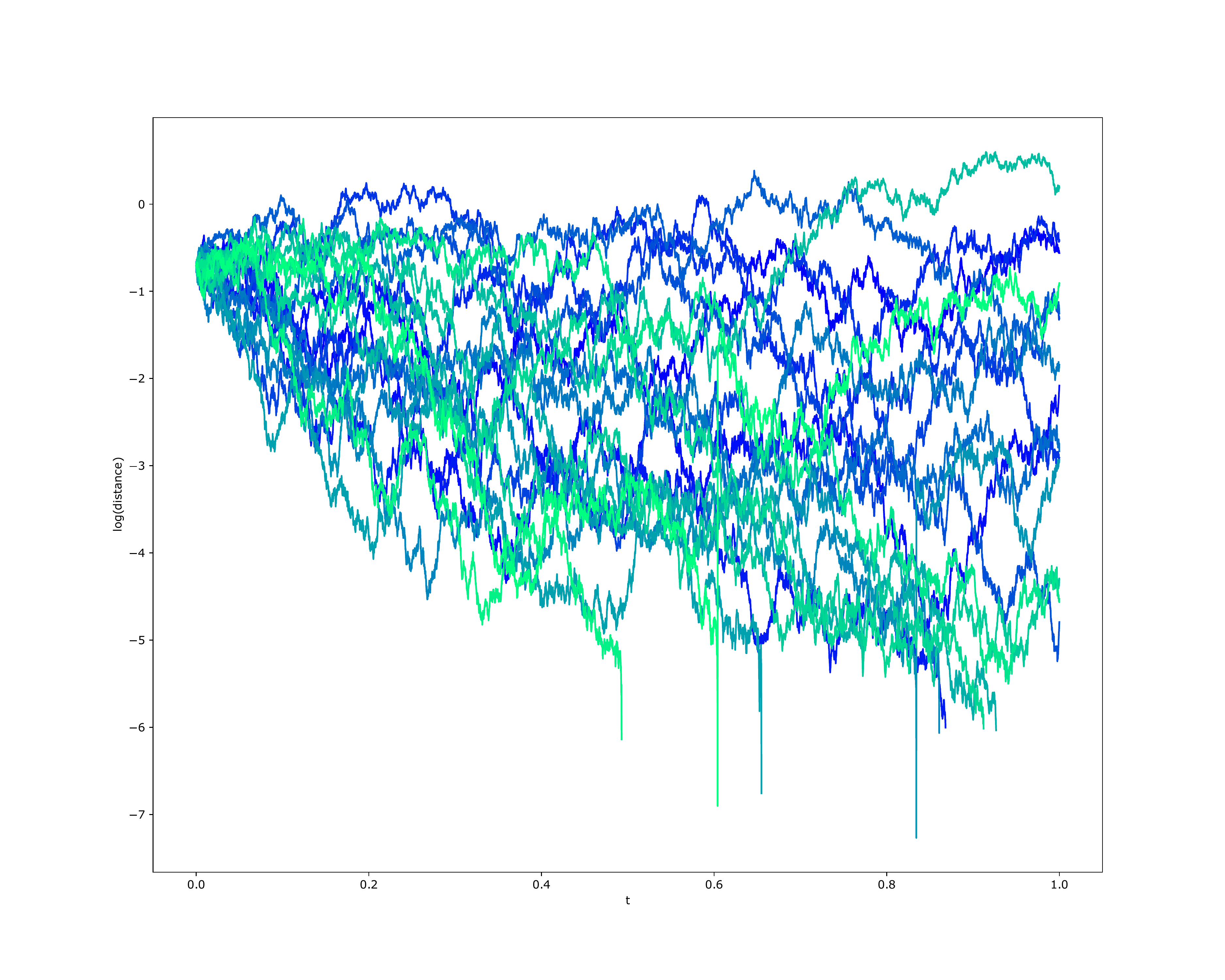}
    \includegraphics[width=0.32\textwidth]{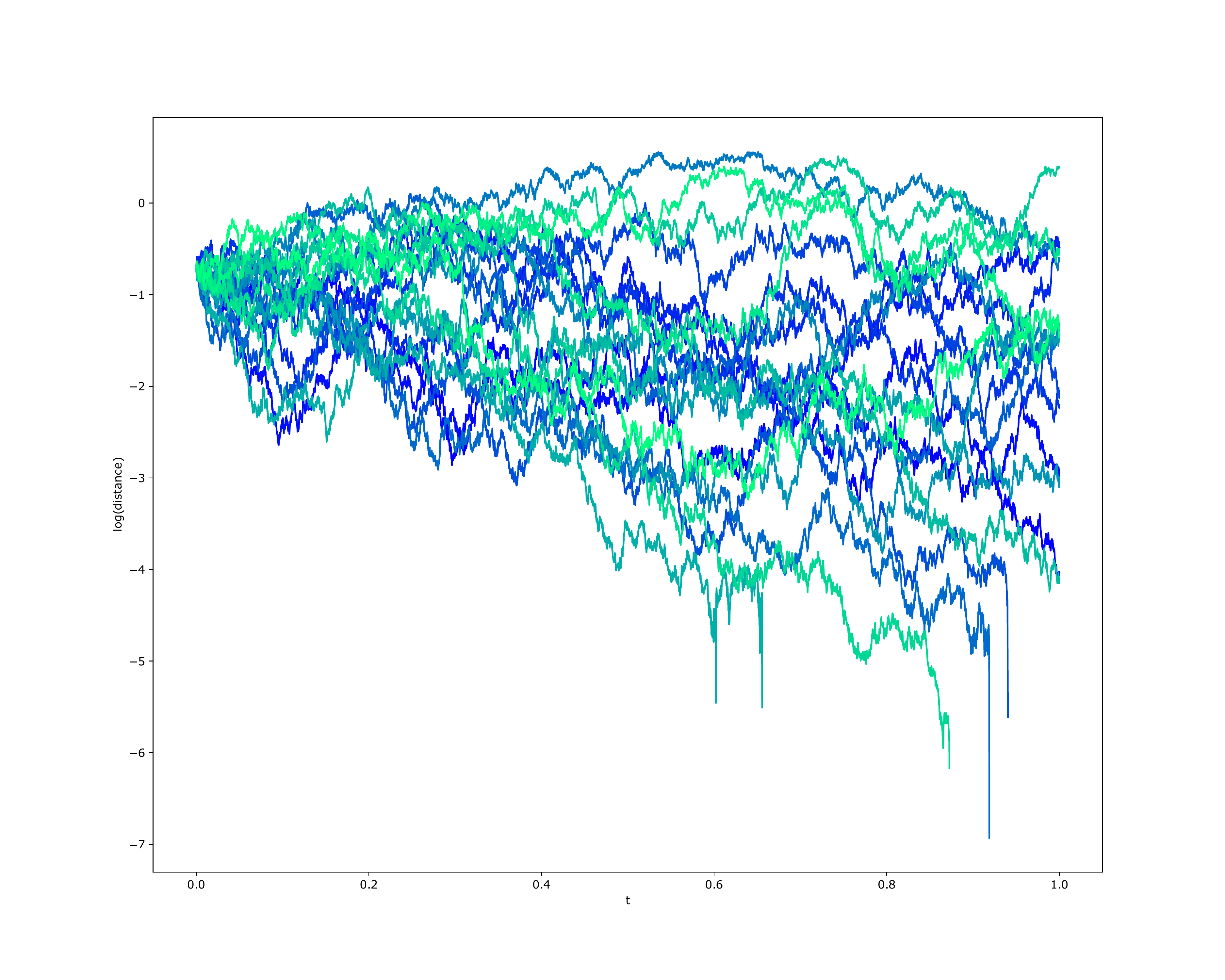}

    \includegraphics[width=0.32\textwidth]{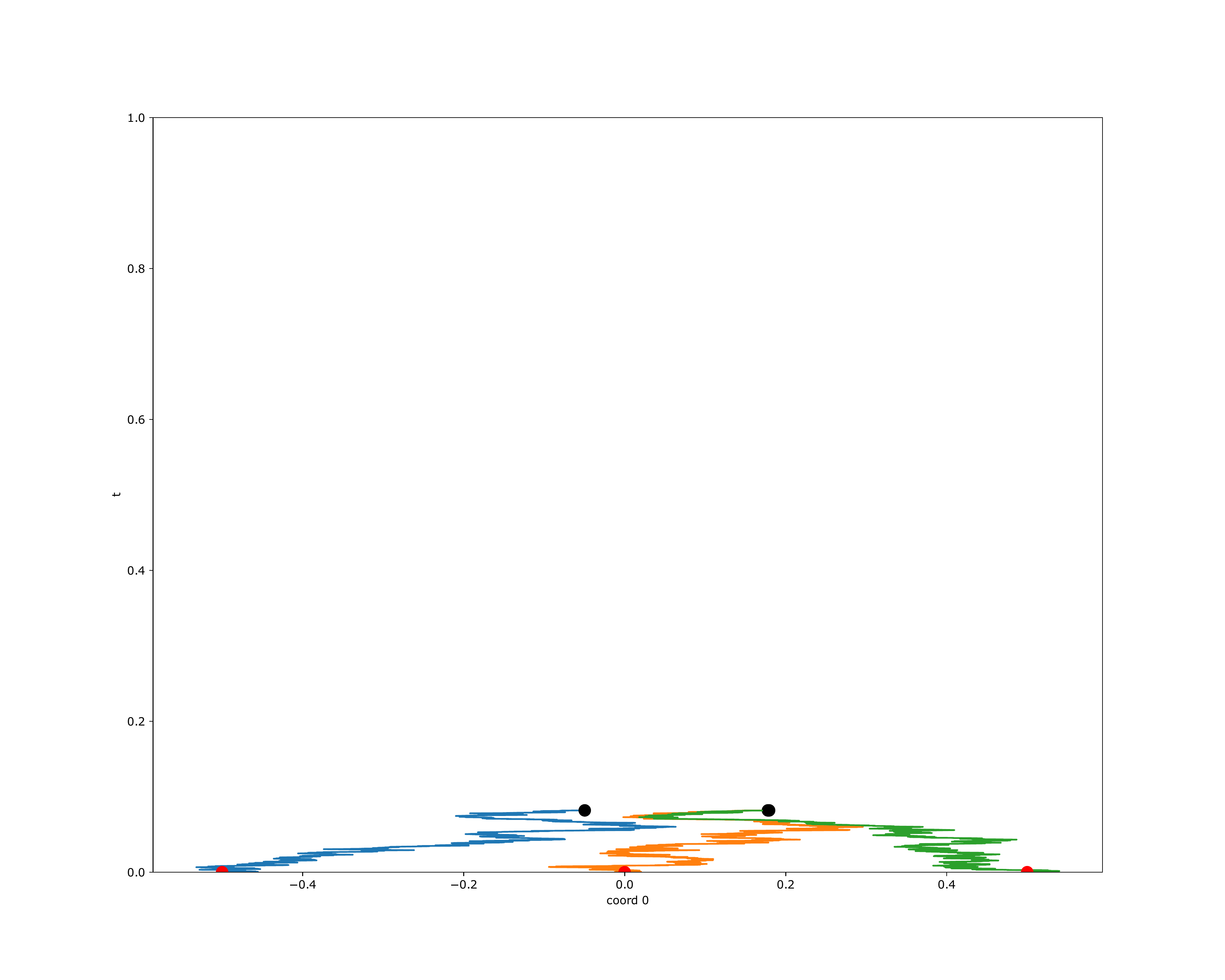}
    \includegraphics[width=0.32\textwidth]{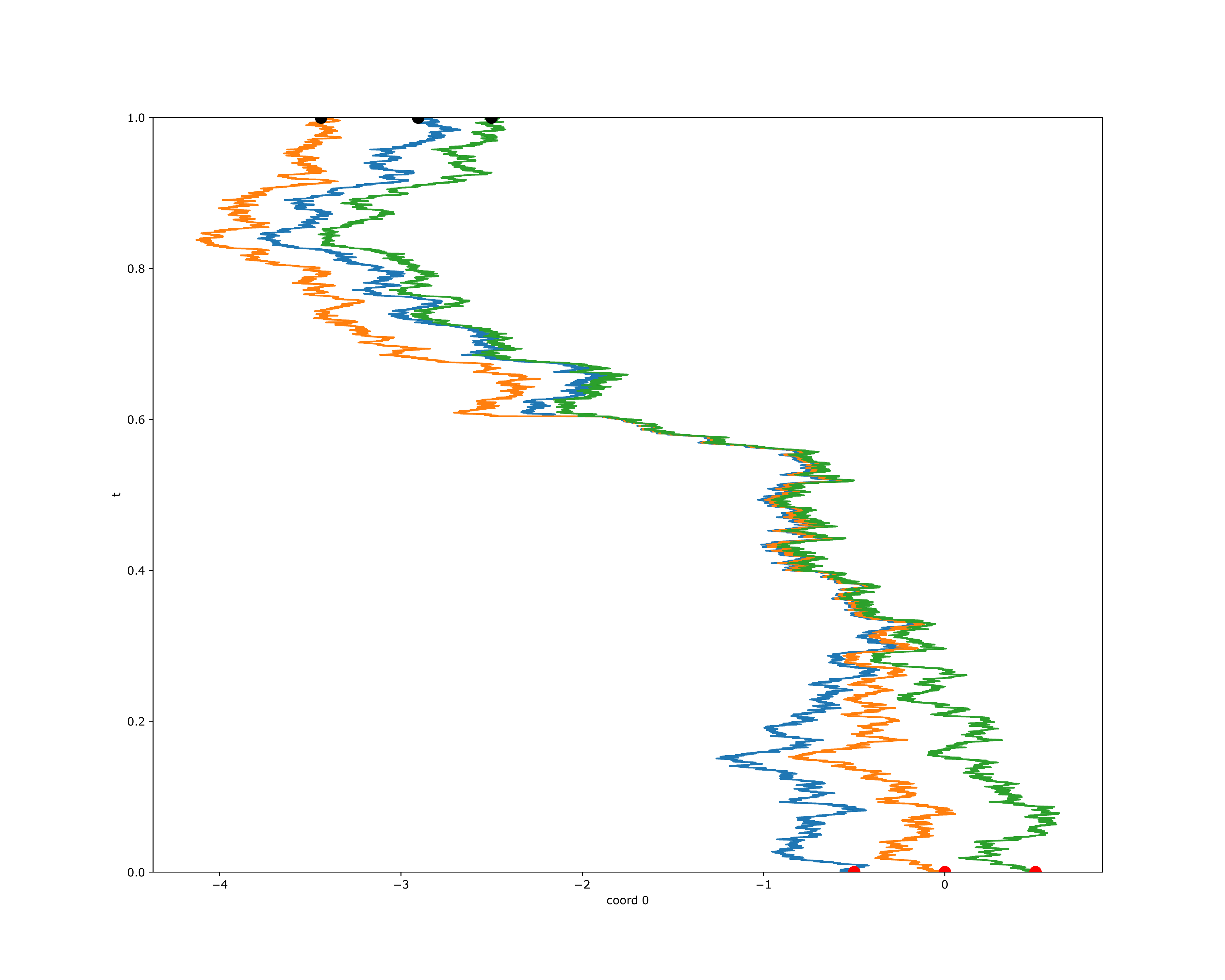}
    \includegraphics[width=0.32\textwidth]{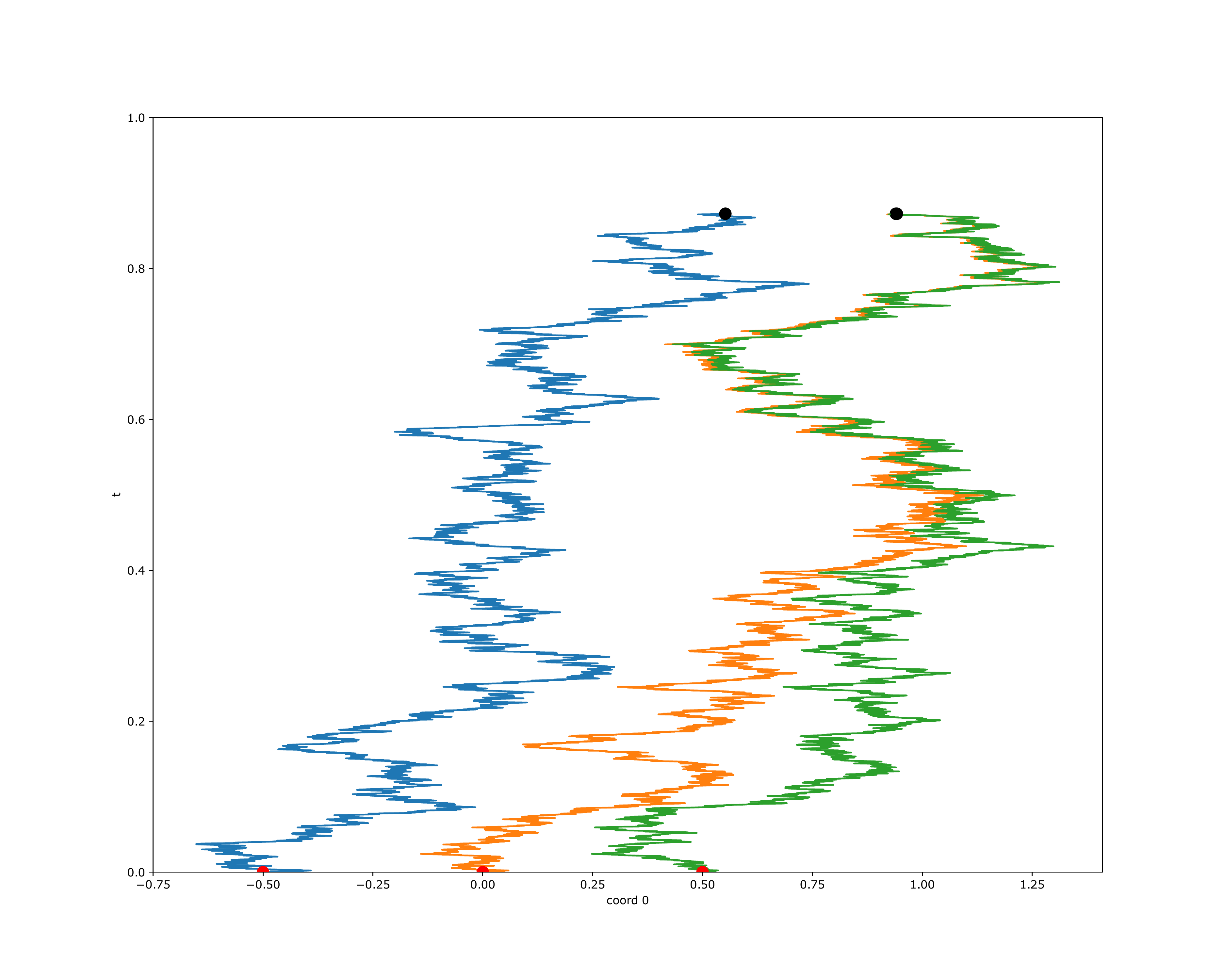}

    \caption{Results for $d=1$ and $n=3$. Setup as in Figure~\ref{fig:1D}.}
    \label{fig:N3}
\end{figure}

\begin{figure}
    \centering
    \includegraphics[width=0.32\textwidth]{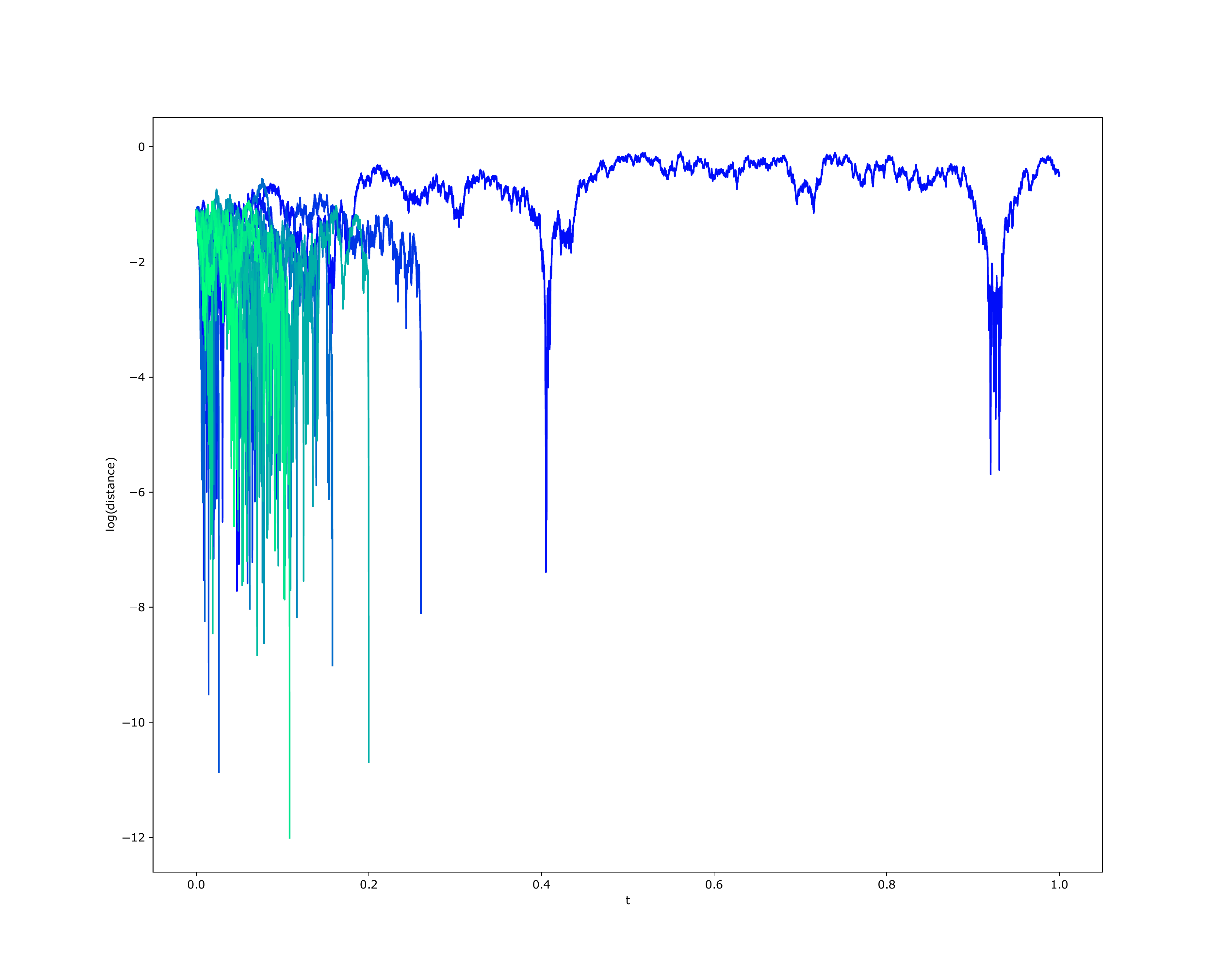}
    \includegraphics[width=0.32\textwidth]{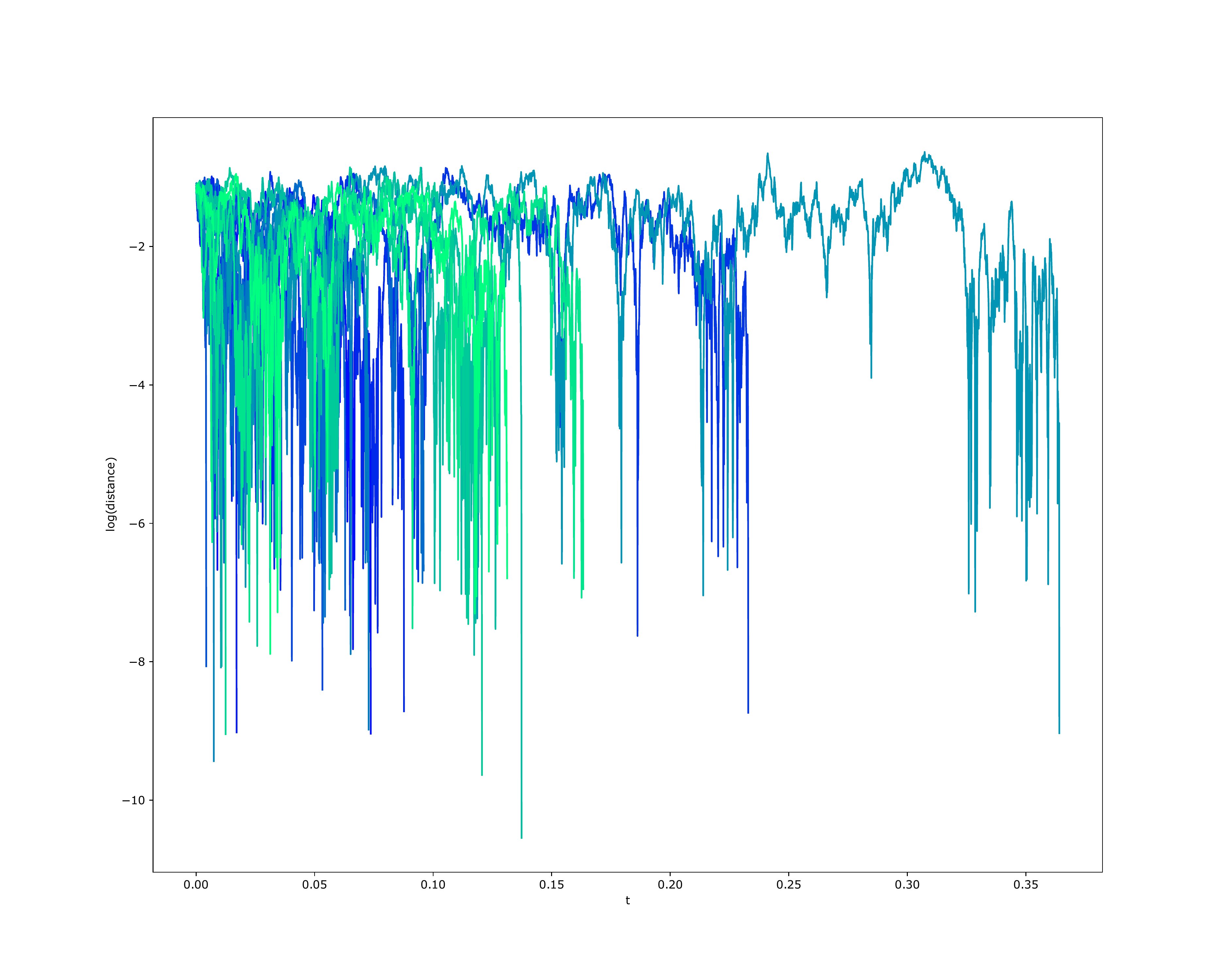}
    \includegraphics[width=0.32\textwidth]{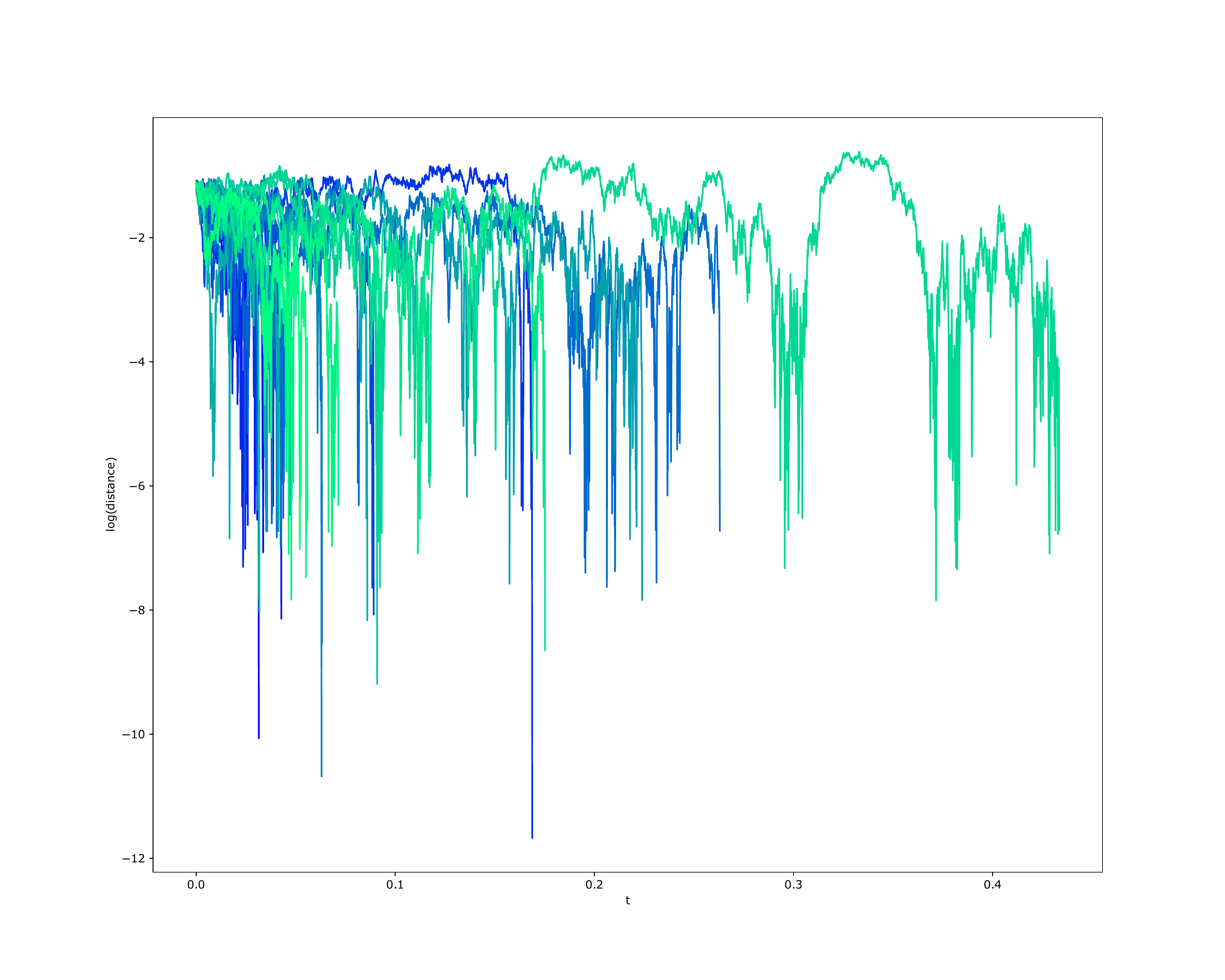}

    \includegraphics[width=0.32\textwidth]{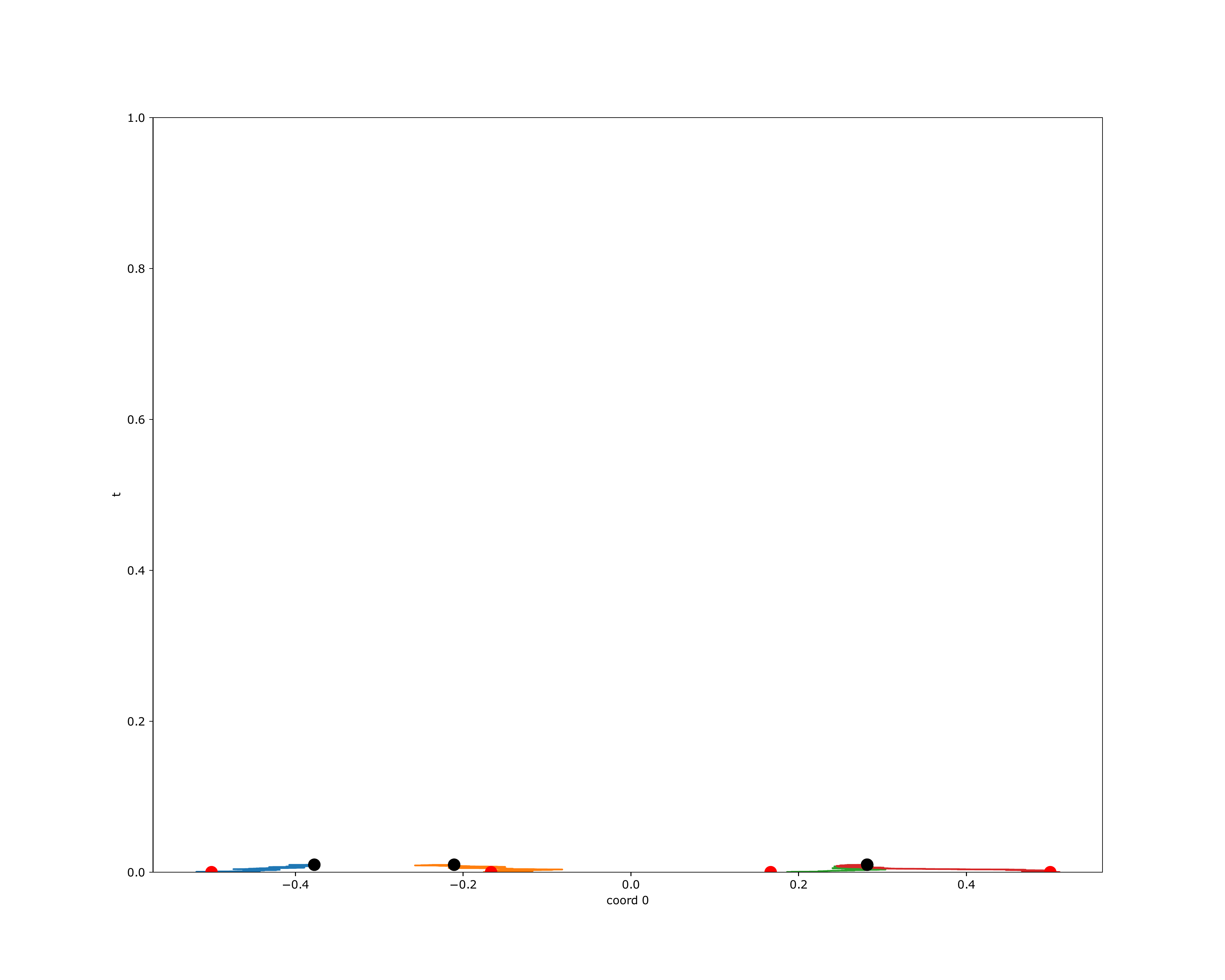}
    \includegraphics[width=0.32\textwidth]{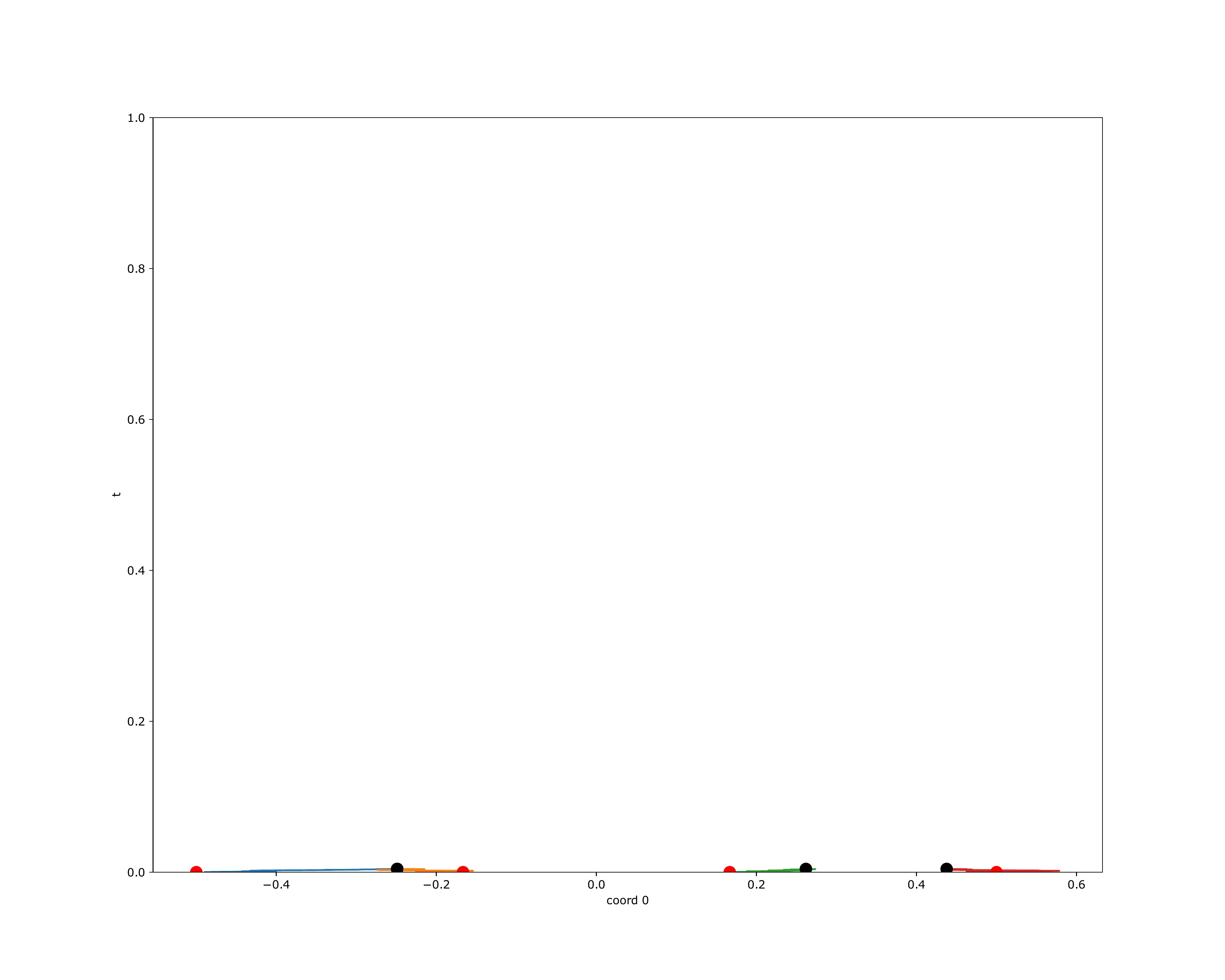}
    \includegraphics[width=0.32\textwidth]{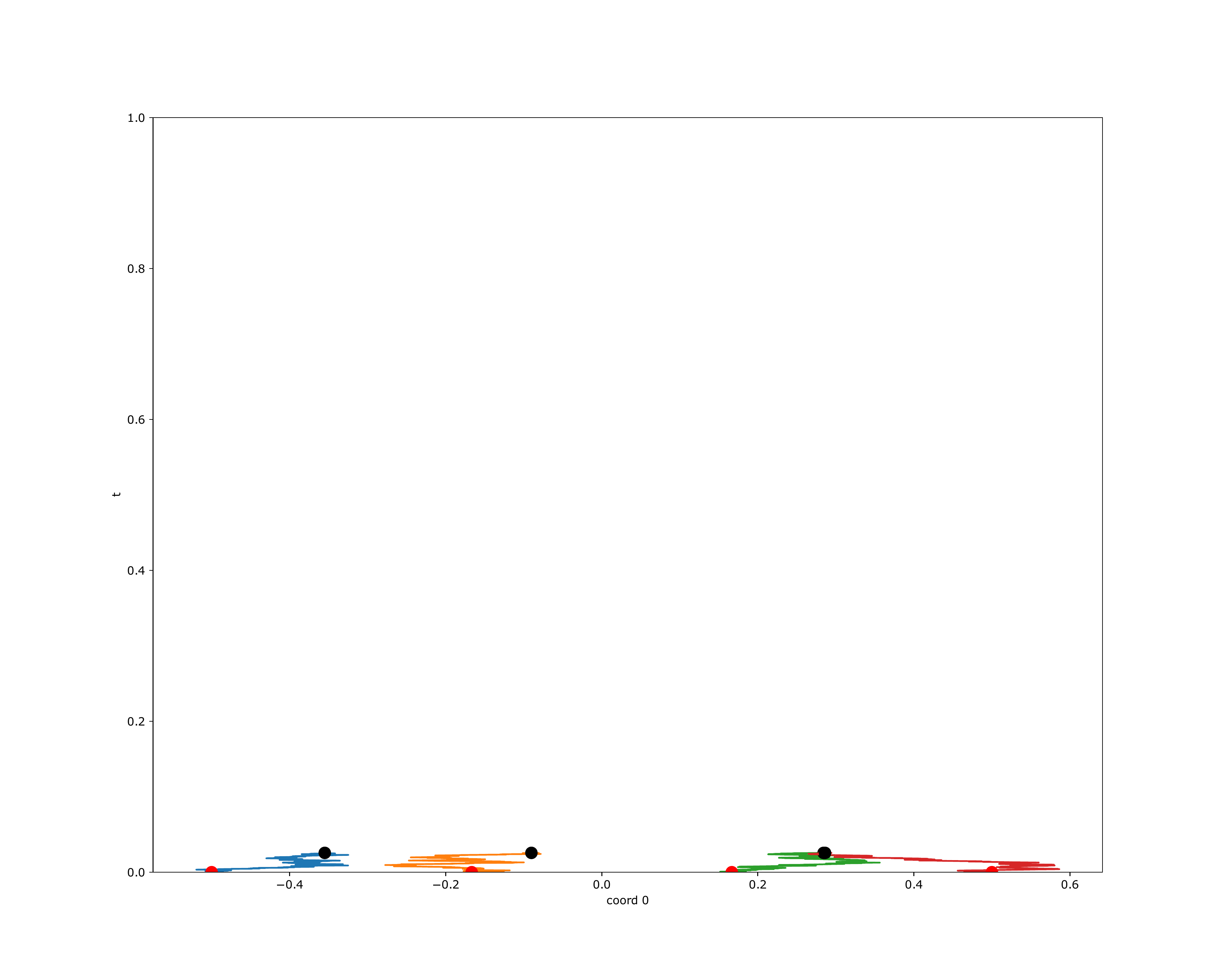}

    \caption{Results for $d=1$ and $n=4$. Setup as in Figure~\ref{fig:1D}.}
    \label{fig:N4}
\end{figure}

\begin{figure}
    \centering
    \includegraphics[width=0.32\textwidth]{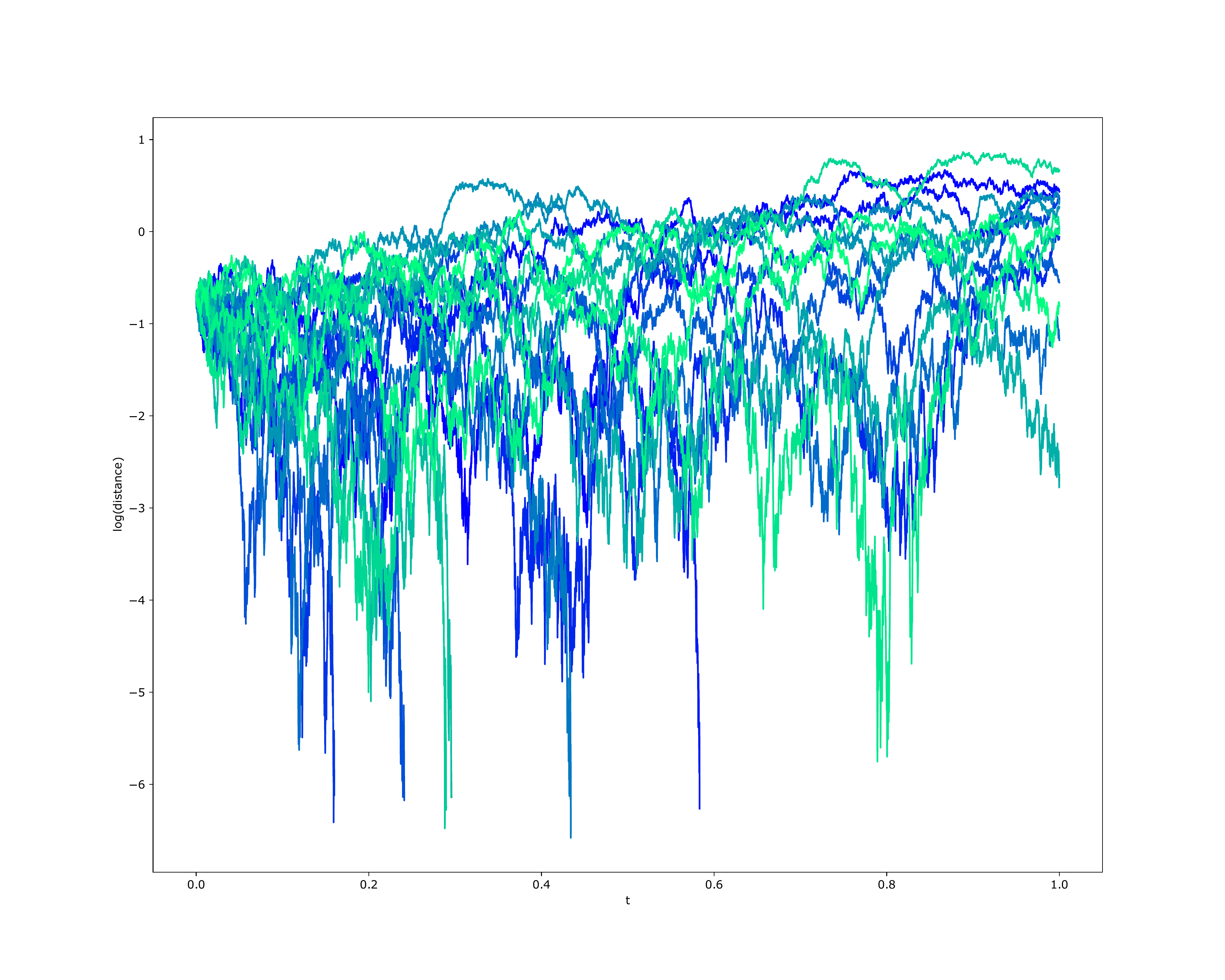}
    \includegraphics[width=0.32\textwidth]{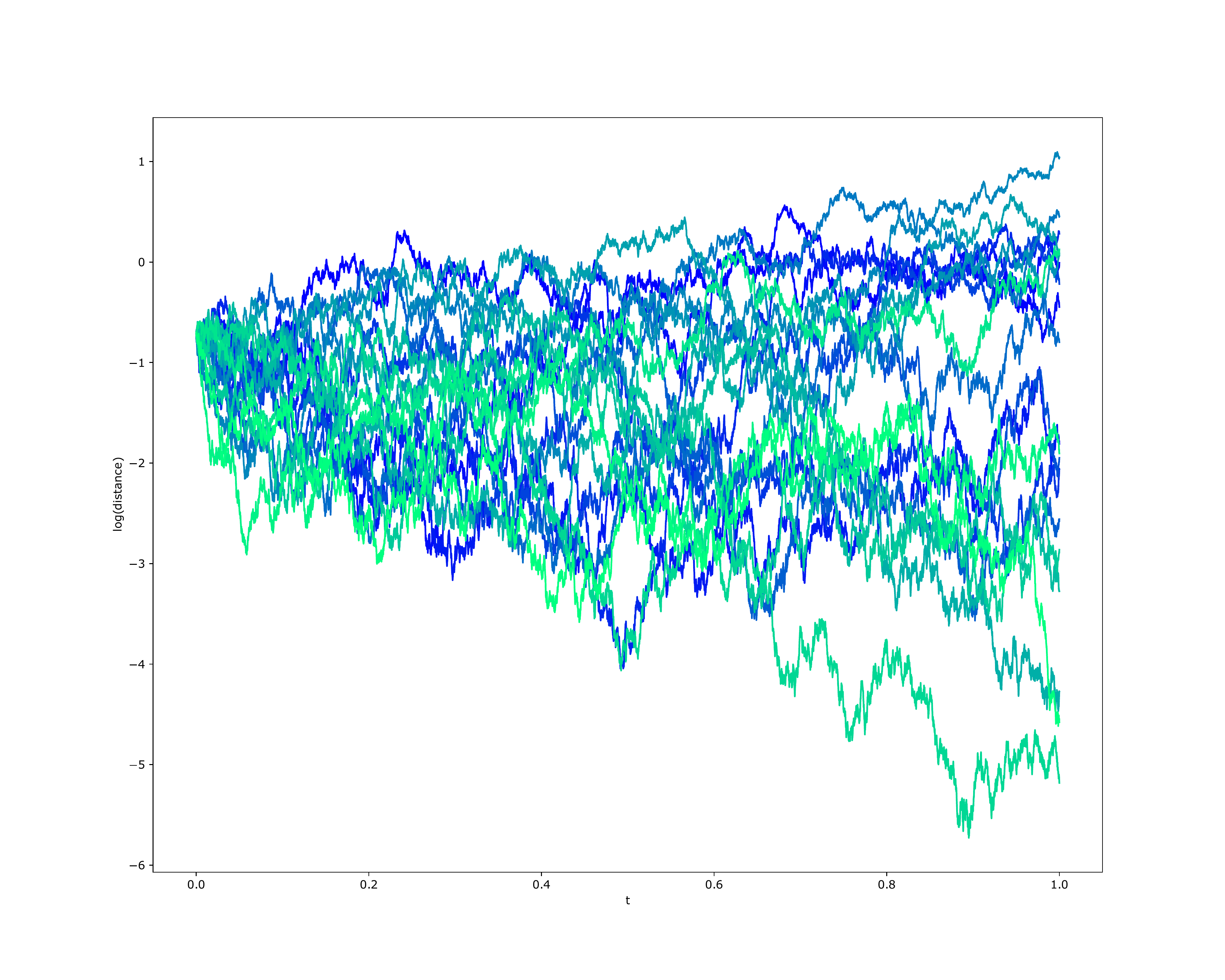}
    \includegraphics[width=0.32\textwidth]{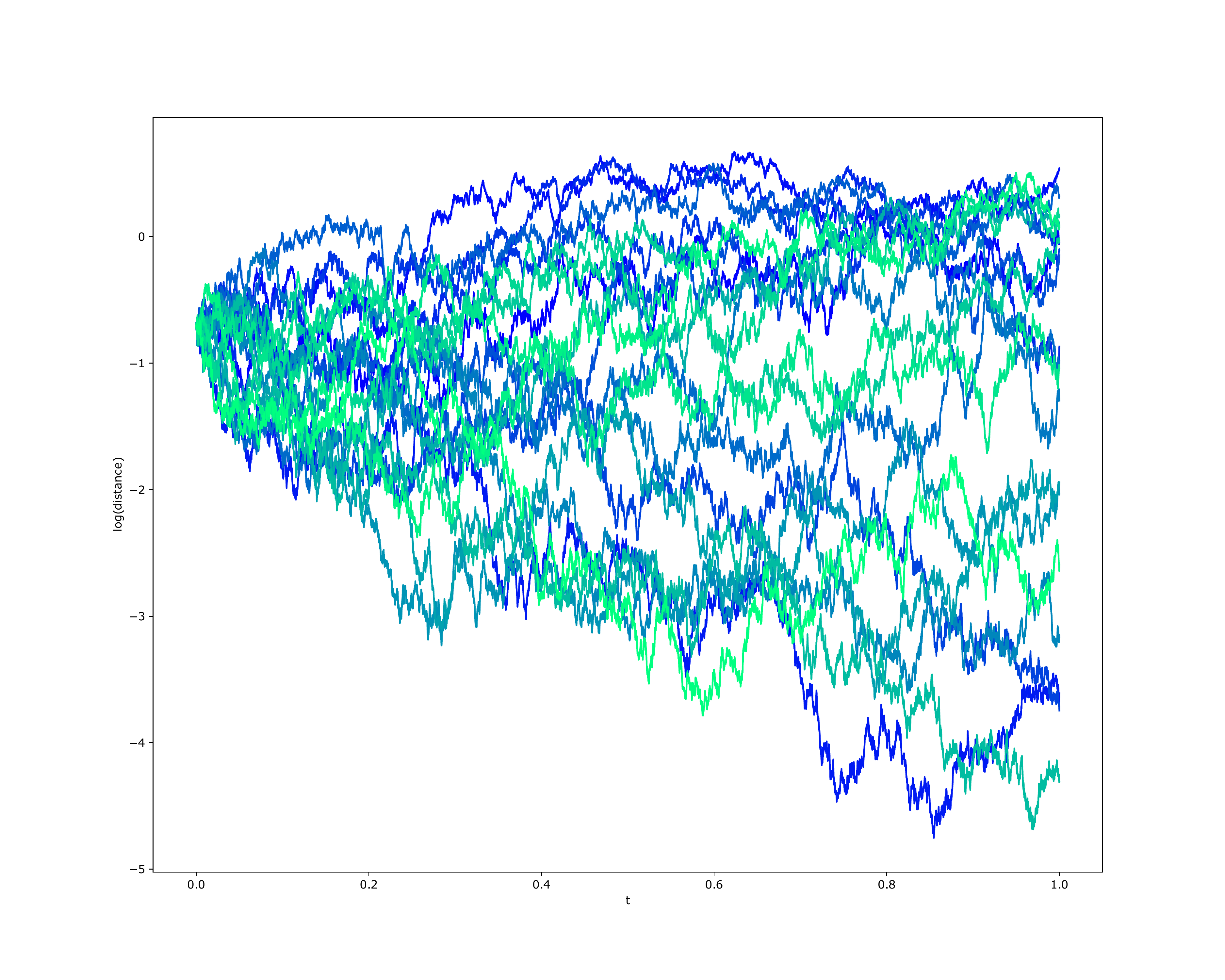}

    \includegraphics[width=0.32\textwidth]{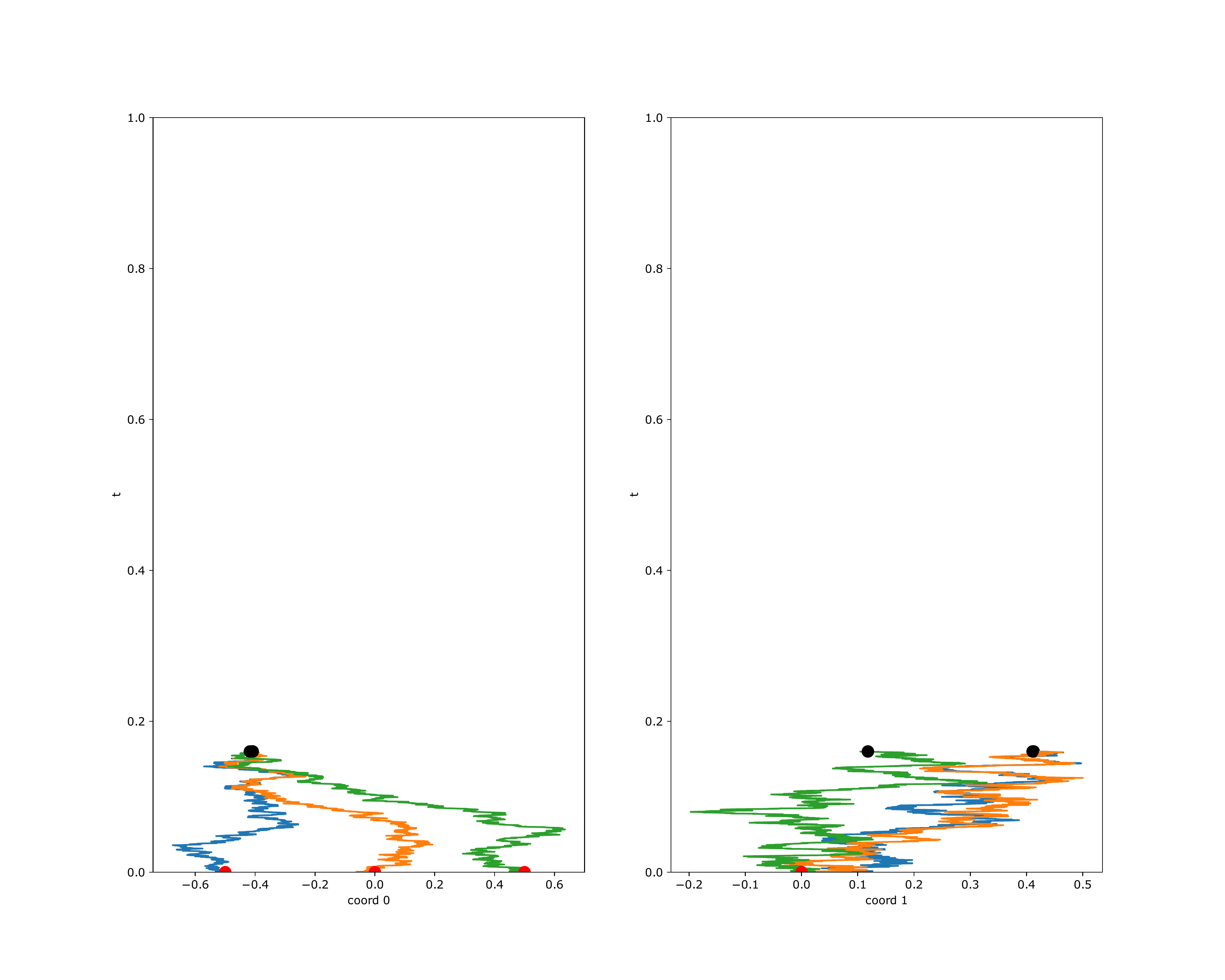}
    \includegraphics[width=0.32\textwidth]{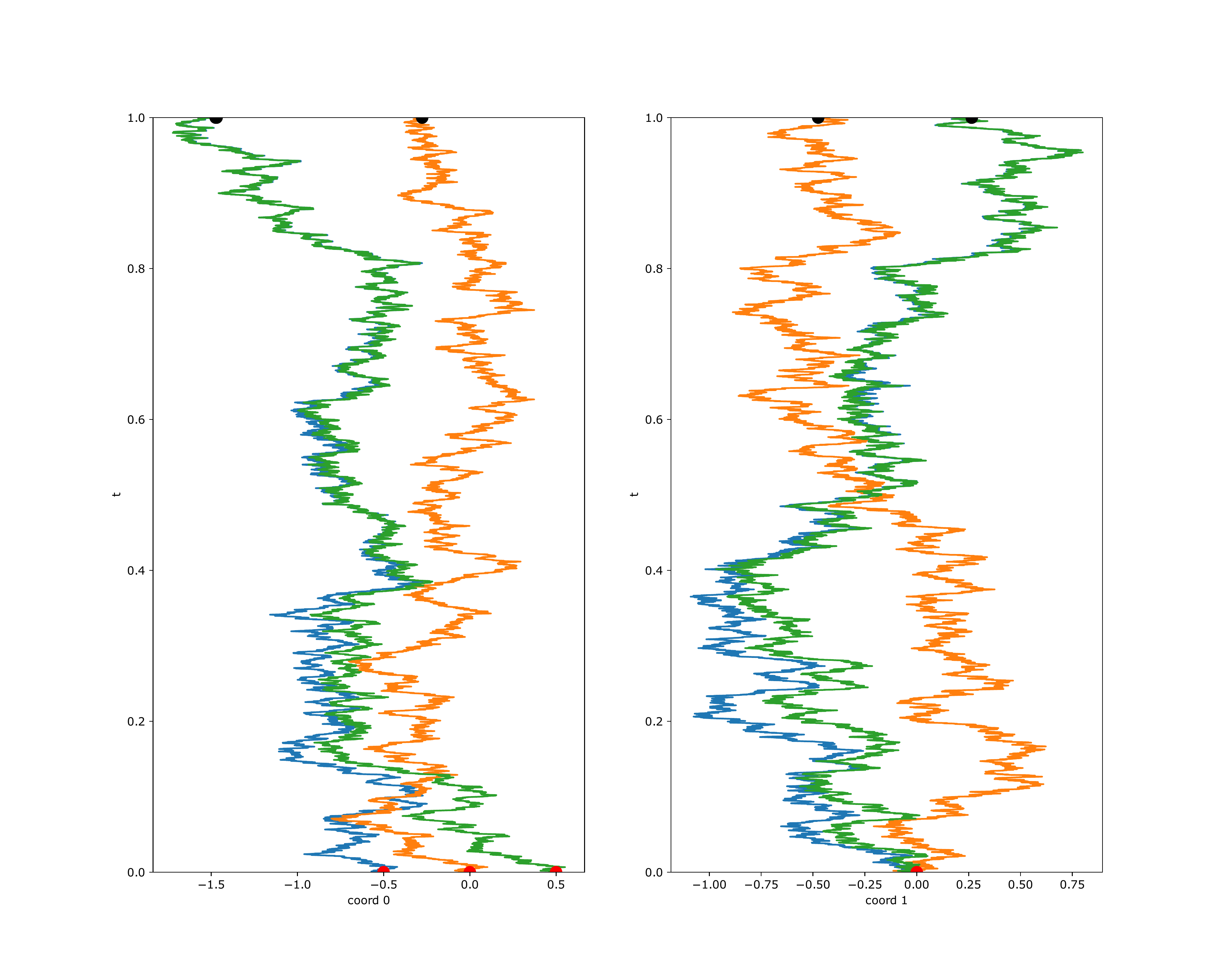}
    \includegraphics[width=0.32\textwidth]{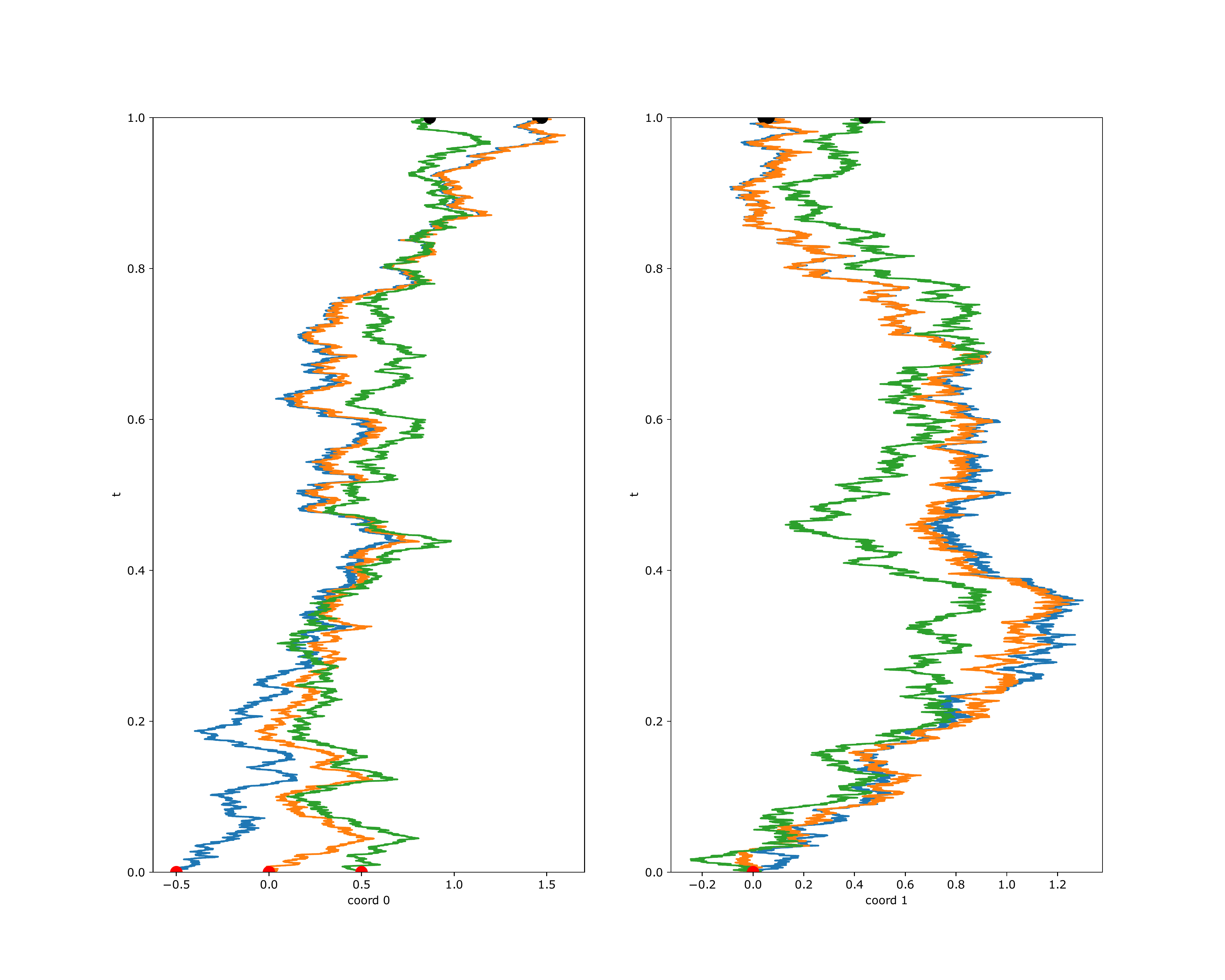}

    \caption{Results for $d=2$ and $n=3$. Setup as in Figure~\ref{fig:1D}.}
    \label{fig:N3D2}
\end{figure}

We now repeat the above experiments for more than two landmarks. 
Figure~\ref{fig:N3} shows $n=3$ landmarks in dimension $d=1$ and indicates a potential collision for the kernel $k_{3/2}$ and the Gaussian kernel, in addition to the kernel $k_{1/2}$. Figure~\ref{fig:N4} shows $n=4$ landmarks in dimension $d=1$. Here, interestingly, collision seems to occur almost immediately for all kernels.
Figure~\ref{fig:N3D2} shows $n=3$ landmarks in dimension $d=2$. The results in this case appear to differ from dimension $d=1$, with collision observed for the kernel $k_{1/2}$ only. It thus seems that dimension plays an important role, for fixed small-distance asymptotics of the kernel.

It must be stressed once again that the numerical experiments are only indicative and not conclusive. Particular care must be taken when interpreting numerical results where distances between landmarks get close to the machine precision.

\clearpage 

\bibliographystyle{plain}
\bibliography{references,library}

\begin{thebibliography}{10}

\bibitem{abramowitz1972handbook}
Milton Abramowitz and Irene~A Stegun.
\newblock {\em Handbook of Mathematical Functions with Formulas, Graphs, and
  Mathematical Tables. National Bureau of Standards Applied Mathematics Series
  55. Tenth Printing.}
\newblock ERIC, 1972.

\bibitem{aronszajn1961theory}
Nachman Aronszajn and Kennan~T Smith.
\newblock Theory of {B}essel potentials. {I}.
\newblock {\em Annales de l'Institut Fourier}, 11:385--475, 1961.

\bibitem{bauer2022smooth}
Martin Bauer, Martins Bruveris, Philipp Harms, and Peter~W Michor.
\newblock Smooth perturbations of the functional calculus and applications to
  {R}iemannian geometry on spaces of metrics.
\newblock {\em Communications in Mathematical Physics}, 389(2):899--931, 2022.

\bibitem{bauerOverviewGeometriesShape2014}
Martin Bauer, Martins Bruveris, and Peter~W Michor.
\newblock Overview of the {{Geometries}} of {{Shape Spaces}} and
  {{Diffeomorphism Groups}}.
\newblock {\em J. Math. Imaging Vis.}, 50(1-2):60--97, September 2014.

\bibitem{bauer2023regularity}
Martin Bauer, Philipp Harms, and Peter~W Michor.
\newblock Regularity and completeness of half-{L}ie groups.
\newblock {\em \textnormal{\texttt{arXiv:2302.01631}}}, 2023.

\bibitem{SSDE}
Alexander Cherny and Hans-J\"{u}rgen Engelbert.
\newblock {\em Singular stochastic differential equations}, volume 1858 of {\em
  Lecture Notes in Mathematics}.
\newblock Springer, Berlin, 2005.

\bibitem{eltznerDiffusionMeansGeometric2022}
Benjamin Eltzner, Pernille Hansen, Stephan~F Huckemann, and Stefan Sommer.
\newblock Diffusion means in geometric spaces.
\newblock {\em \textnormal{\texttt{arXiv:2105.12061}}}, 2022.

\bibitem{hackenbroch1994stochastische}
Wolfgang Hackenbroch and Anton Thalmaier.
\newblock {\em Stochastische Analysis}.
\newblock Teubner, Stuttgart, 1994.

\bibitem{hsu}
Elton~P Hsu.
\newblock {\em Stochastic {A}nalysis on {M}anifolds}, volume~38 of {\em
  Graduate Studies in Mathematics}.
\newblock American Mathematical Society, Providence, RI, 2002.

\bibitem{joshi2000landmark}
Sarang~C Joshi and Michael~I Miller.
\newblock Landmark matching via large deformation diffeomorphisms.
\newblock {\em IEEE transactions on image processing}, 9(8):1357--1370, 2000.

\bibitem{kriegl2015exotic}
Andreas Kriegl, Peter~W Michor, and Armin Rainer.
\newblock An exotic zoo of diffeomorphism groups on {$\mathbb R^n$}.
\newblock {\em Annals of Global Analysis and Geometry}, 47:179--222, 2015.

\bibitem{micheliDifferentialGeometryLandmark2008}
Mario Micheli.
\newblock {\em The differential geometry of landmark shape manifolds: metrics,
  geodesics, and curvature}.
\newblock PhD thesis, Brown University, Providence, RI, 2008.

\bibitem{micheli2012sectional}
Mario Micheli, Peter~W Michor, and David Mumford.
\newblock Sectional curvature in terms of the cometric, with applications to
  the {R}iemannian manifolds of landmarks.
\newblock {\em SIAM Journal on Imaging Sciences}, 5(1):394--433, 2012.

\bibitem{michor2020manifolds}
Peter~W Michor.
\newblock {\em Geometric Continuum Mechanics}, volume~42 of {\em Advances in
  Continuum Mechanics}, chapter Manifolds of mappings for continuum mechanics,
  pages 3--75.
\newblock Birkh{\"a}user, Basel, 2020.

\bibitem{michor2013zoo}
Peter~W Michor and David Mumford.
\newblock A zoo of diffeomorphism groups on $\mathbb{R}^n$.
\newblock {\em Annals of Global Analysis and Geometry}, 44:529--540, 2013.

\bibitem{sommerBridgeSimulationMetric2017}
Stefan Sommer, Alexis Arnaudon, Line Kuhnel, and Sarang Joshi.
\newblock Bridge simulation and metric estimation on landmark manifolds.
\newblock In {\em Graphs in {{Biomedical Image Analysis}}, {{Computational
  Anatomy}} and {{Imaging Genetics}}}, Lecture {{Notes}} in {{Computer
  Science}}, pages 79--91. {Springer}, September 2017.

\bibitem{stanevaLearningShapeTrends2017}
Valentina Staneva and Laurent Younes.
\newblock Learning {{Shape Trends}}: {{Parameter Estimation}} in {{Diffusions}}
  on {{Shape Manifolds}}.
\newblock In {\em 2017 {{IEEE Conference}} on {{Computer Vision}} and {{Pattern
  Recognition Workshops}} ({{CVPRW}})}, pages 717--725, July 2017.

\bibitem{younes}
Laurent Younes.
\newblock {\em Shapes and diffeomorphisms}, volume 171 of {\em Applied
  Mathematical Sciences}.
\newblock Springer, Berlin, 2010.

\end{thebibliography}

\end{document}